\documentclass[aos,nameyear,dvips]{arximspdf}
\usepackage{graphicx}
\doi{10.1214/09-AOS772}
\volume{38}
\issue{6}
\pubyear{2010}
\firstpage{3412}
\lastpage{3444}

\makeatletter
\newtheorem{theorem}{Theorem}
\newtheorem{proposition}[theorem]{Proposition}
\newtheorem{lemma}[theorem]{Lemma}
\newtheorem{corollary}[theorem]{Corollary}

\newproclaim{example}{Example}
\newproclaim{conjecture}{Conjecture}
\newproclaim{definition}{Definition}
\newproclaim{remark}{Remark}
\newproclaim{assumption}{Assumption}
\newproclaim{condition}{Condition}

\def\afrac#1#2{#1/(#2)}
\def\vafrac#1#2{(#1)/(#2)}

\makeatother

\begin{document}
\begin{frontmatter}

\title{A reproducing kernel Hilbert space approach to functional linear regression}
\runtitle{RKHS approach to FLR}

\begin{aug}
\author[A]{\fnms{Ming} \snm{Yuan}\corref{}\thanksref{t1}\ead[label=e1]{myuan@isye.gatech.edu}} and
\author[B]{\fnms{T. Tony} \snm{Cai}\thanksref{t2}\ead[label=e2]{tcai@wharton.upenn.edu}}
\thankstext{t1}{Supported in part by NSF Grant DMS-MPSA-0624841 and NSF
CAREER Award DMS-0846234.}
\thankstext{t2}{Supported in part by NSF Grant DMS-0604954 and NSF FRG
Grant DMS-0854973.}
\runauthor{M. Yuan and T. T. Cai}
\affiliation{Georgia Institute of Technology and University of Pennsylvania}
\address[A]{Milton Stewart School of Industrial and\\
\quad Systems Engineering\\
Georgia Institute of Technology\\
Atlanta, Georgia 30332\\
USA\\
\printead{e1}} %adresu isvedimo komanda gale!
\address[B]{Department of Statistics\\
The Wharton School\\
University of Pennsylvania\\
Philadelphia, Pennsylvania 19104\\
USA\\
\printead{e2}}
\end{aug}

% HISTORY:
\received{\smonth{3} \syear{2009}}
\revised{\smonth{11} \syear{2009}}

% ABSTRACT
%
\begin{abstract}
We study in this paper a smoothness regularization method for
functional linear regression and provide a unified treatment for both
the prediction and estimation problems. By
developing a tool on simultaneous diagonalization of two positive
definite kernels, we obtain shaper results on the minimax rates of
convergence and show that smoothness regularized estimators
achieve the optimal rates of convergence for both prediction and
estimation under conditions weaker than those for the
functional principal components based methods developed in the literature.
Despite the generality of the method of regularization, we show that
the procedure is easily implementable. Numerical results are obtained to
illustrate the merits of the method and to demonstrate the theoretical
developments.
\end{abstract}

% KEYWORDS
%

\end{frontmatter}

%s1 ###
\section{Introduction}
Consider the following functional linear regression mod\-el where the
response $Y$ is related to a square integrable random function $X(\cdot
)$ through
%
%e1 ###
%
\begin{equation}
Y=\alpha_0+\int_\mathcal{T}X(t)\beta_0(t)\,dt+\varepsilon.
\end{equation}
Here $\alpha_0$ is the intercept, $\mathcal{T}$ is the domain of
$X(\cdot)$,
$\beta_0(\cdot)$ is an unknown slope function and $\varepsilon$ is
a centered noise random variable. The domain
$\mathcal{T}$ is assumed to be a compact subset of an Euclidean space.
Our goal is to estimate
$\alpha_0$ and $\beta_0(\cdot)$ as well as to retrieve
%
%e2 ###
%
\begin{equation}
\label{eq:lr}
\eta_0(X):=\alpha_0+\int_\mathcal{T}X(t)\beta_0(t)\,dt
\end{equation}
based on a set of training data $(x_1,y_1),\ldots,(x_n,y_n)$
consisting of $n$ independent copies of $(X,Y)$. We shall assume that
the slope function $\beta_0$ resides in a reproducing kernel Hilbert
space (RKHS) $\mathcal{H}$, a subspace of the collection of square integrable
functions on $\mathcal{T}$.

In this paper, we investigate the method of regularization for
estimating $\eta_0$, as well as $\alpha_0$ and $\beta_0$. Let $\ell_n$
be a data fit functional that measures how well $\eta$ fits the data
and $J$ be a penalty functional that assesses the ``plausibility'' of
$\eta$. The method of regularization estimates $\eta_0$ by
%
%e3 ###
%
\begin{equation}
\label{eq:mor}
\hat{\eta}_{n\lambda}=\mathop{\arg\min}_{\eta} [\ell_n (\eta
|{\rm
data} )+\lambda J(\eta) ],
\end{equation}
where the minimization is taken over
%
%e4 ###
%
\begin{equation}
\biggl\{\eta\dvtx\mathcal{L}_2(\mathcal{T}) \to\mathbb{R}\big| \eta
(X)=\alpha+\int
_\mathcal{T}X\beta\dvtx
\alpha\in\mathbb{R}, \beta\in\mathcal{H}\biggr\},
\end{equation}
and $\lambda\ge0$ is a tuning parameter that balances the fidelity to
the data and the plausibility. Equivalently, the minimization can be
taken over $(\alpha,\beta)$ instead of $\eta$ to obtain estimates for
both the intercept and slope, denoted by $\hat{\alpha}_{n\lambda}$ and
$\hat{\beta}_{n\lambda}$ hereafter. The most common choice of the data
fit functional is the squared error
%
%e5 ###
%
\begin{equation}
\label{eq:ls}
\ell_n(\eta)={1\over n}\sum_{i=1}^n [y_i-\eta(x_i) ]^2.
\end{equation}
In general, $\ell_n$ is chosen such that it is convex in $\eta$ and
$E\ell_n(\eta)$ in uniquely minimized by $\eta_0$.

In the context of functional linear regression, the penalty functional
can be conveniently defined through the slope function $\beta$ as a
squared norm or semi-norm associated with $\mathcal{H}$. The canonical
example of
$\mathcal{H}$ is the Sobolev spaces. Without loss of generality,
assume that
$\mathcal{T}=[0,1]$, the Sobolev space of order $m$ is then defined as
\begin{eqnarray*}
\mathcal{W}_2^m([0,1])&=& \bigl\{\beta\dvtx[0,1]\to\mathbb{R}|
\beta,
\beta^{(1)},\ldots,\beta^{(m-1)} \mbox{ are absolutely}
\\
&&\hspace*{104pt}{}\mbox{continuous and } \beta^{(m)}\in\mathcal{L}_2
\bigr\}.
\end{eqnarray*}
There are many possible norms that can be equipped with $\mathcal
{W}_2^m$ to
make it a reproducing kernel Hilbert space. For example, it can be
endowed with the norm
%
%e6 ###
%
\begin{equation}
\label{eq:sobnorm}
\Vert\beta\Vert^2_{\mathcal{W}_2^m}=\sum_{q=0}^{m-1} \biggl(\int
\beta
^{(q)} \biggr)^2+\int\bigl(\beta^{(m)} \bigr)^2.
\end{equation}
The readers are referred to Adams (\citeyear{Adams1975}) for a
thorough treatment of
this subject. In this case, a possible\vadjust{\goodbreak} choice of the penalty functional
is given~by\looseness=1
%
%e7 ###
%
\begin{equation}
\label{eq:sobpen}
J(\beta)=\int_0^1 \bigl[\beta^{(m)}(t) \bigr]^2\,dt.
\end{equation}\looseness=0
Another setting of particular interest is $\mathcal{T}=[0,1]^2$ which
naturally occurs when $X$ represents an image. A popular choice
in this setting is the thin plate spline where $J$ is given by
%
%e8 ###
%
\begin{equation}
J(\beta)=\int_0^1\!\!\int_0^1
\biggl[ \biggl({\partial^2 \beta\over\partial x_1^2} \biggr)^2
+2 \biggl({\partial^2 \beta\over\partial x_1\,\partial x_2} \biggr)^2
+ \biggl({\partial^2 \beta\over\partial x_2^2} \biggr)^2 \biggr]\,
dx_1\,dx_2,
\end{equation}
and $(x_1,x_2)$ are the arguments of bivariate function $\beta$. Other
examples of $\mathcal{T}$ include $\mathcal{T}=\{1,2,\ldots, p\}$
for some
positive integer $p$, and unit sphere in an Euclidean space among
others. The readers are referred to Wahba (\citeyear{Wahba1990}) for
common choices of
$\mathcal{H}$ and $J$ in these as well as other contexts.

Other than the methods of regularization, a number of alternative
estimators have been introduced in recent years for the functional
linear regression [James (\citeyear{James2002});
Cardot, Ferraty and Sarda (\citeyear{CardotFerratySarda2003});
Ramsay and Silverman (\citeyear{RamsaySilverman2005});
Yao, M\"{u}ller and Wang (\citeyear{YaoMullerWang2005});
Ferraty and Vieu (\citeyear{FerratyVieu2006});
Cai and Hall (\citeyear{CaiHall2006});
Li and Hsing (\citeyear{LiHsing2007});
Hall and Horowitz (\citeyear{HallHorowitz2007});
Crambes, Kneip and Sarda (\citeyear{CrambesKneipSarda2009});
Johannes (\citeyear{Johanness2009})]. Most
of the existing
methods are based upon the functional principal component analysis
(FPCA). The success of these approaches hinges on the availability of a
good estimate of the functional principal components for $X(\cdot)$. In
contrast, the aforementioned smoothness regularized estimator avoids
this task and therefore circumvents assumptions on the spacing of the
eigenvalues of the covariance operator for $X(\cdot)$ as well as
Fourier coefficients of $\beta_0$ with respect to the eigenfunctions,
which are required by the FPCA-based approaches. Furthermore, as we
shall see in the subsequent theoretical analysis, because the
regularized estimator does not rely on estimating the functional
principle components, stronger results on the convergence rates can be obtained.

Despite the generality of the method of regularization, we show that
the estimators can be computed rather efficiently.
We first derive a representer
theorem in Section~\ref{representer.sec} which demonstrates that
although the minimization with respect to $\eta$ in (\ref{eq:mor}) is
taken over an infinite-dimensional space, the solution can actually be
found in a finite-dimensional subspace. This result makes our
procedure easily implementable and enables us to take advantage of the
existing techniques and algorithms for smoothing splines to compute
$\hat{\eta}_{n\lambda}$, $\hat{\beta}_{n\lambda}$ and
$\hat\alpha_{n\lambda}$.

We then consider in Section~\ref{diagonal.sec} the relationship
between the eigen structures of the covariance operator for $X(\cdot)$
and the reproducing kernel of the RKHS $\mathcal{H}$. These eigen structures
play prominent roles in determining the difficulty of the prediction
and estimation problems in functional
linear regression. We prove in Section~\ref{diagonal.sec} a result on
simultaneous diagonalization of the reproducing\vadjust{\goodbreak} kernel
of the RKHS $\mathcal{H}$ and the covariance operator of $X(\cdot)$ which
provides a powerful
machinery for studying the minimax rates of convergence.

Section~\ref{rate.sec} investigates the rates of convergence of the
smoothness regularized estimators. Both the minimax upper and lower bounds
are established. The optimal convergence rates are derived in terms
of a class of intermediate norms which provide a wide range of
measures for the estimation accuracy. In particular, this approach gives
a unified treatment for both the prediction of $\eta_0(X)$ and the
estimation of $\beta_0$. The results show that the smoothness
regularized estimators achieve the optimal rate of convergence for
both prediction and estimation under conditions weaker than those for
the functional principal components based methods developed in the
literature.

The representer theorem makes the regularized estimators easy to implement.
Several efficient algorithms are available in the literature that can
be used for the numerical implementation of our procedure.
Section~\ref{numerical.sec}
presents numerical studies to illustrate the merits of the method as
well as demonstrate the theoretical developments. All proofs are
relegated to
Section~\ref{proof.sec}.\vspace*{-3pt}

%s2 ###
\section{Representer theorem}
\label{representer.sec}

The smoothness regularized estimators $\hat{\eta}_{n\lambda}$ and
$\hat\beta_{n\lambda}$ are defined as the solution to a minimization
problem over an infinite-dimensional space. Before studying the
properties of the estimators, we first show that the minimization is
indeed well defined and easily computable thanks to a version of the
so-called representer theorem.

Let the penalty functional $J$ be a squared semi-norm on $\mathcal{H}$ such
that the null space
%
%e9 ###
%
\begin{equation}
\mathcal{H}_0:= \{\beta\in\mathcal{H}\dvtx J(\beta)=0 \}
\end{equation}
is a finite-dimensional linear subspace of $\mathcal{H}$ with orthonormal
basis $\{\xi_1,\ldots,\break \xi_N\}$ where $N:=\operatorname
{dim}(\mathcal{H}_0)$. Denote by
$\mathcal{H}_1$ its orthogonal complement in $\mathcal{H}$ such that
$\mathcal{H}=\mathcal{H}_0\oplus\mathcal{H}_1$. Similarly, for any function
$f\in\mathcal{H}$, there exists a unique decomposition $f=f_0+f_1$
such that
$f_0\in\mathcal{H}_0$ and $f_1\in\mathcal{H}_1$. Note $\mathcal
{H}_1$ forms a
reproducing kernel Hilbert space with the inner product of $\mathcal{H}$
restricted to $\mathcal{H}_1$. Let $K(\cdot,\cdot)$ be the corresponding
reproducing kernel of $\mathcal{H}_1$ such that $J(f_1)=\|f_1\|
_{K}^2=\|f_1\|
_\mathcal{H}^2$
for any $f_1\in\mathcal{H}_1$. Hereafter we use the subscript $K$ to
emphasize the correspondence between the inner product and its
reproducing kernel.

In what follows, we shall assume that $K$ is continuous and square
integrable. Note that $K$ is also a nonnegative definite operator on
$\mathcal{L}_2$. With slight abuse of notation, write
%
%e10 ###
%
\begin{equation}
(Kf)(\cdot)=\int_\mathcal{T}K(\cdot,s)f(s)\,ds.
\end{equation}
It is known [see, e.g., Cucker and Smale (\citeyear{CuckerSmale2001})]
that $Kf\in\mathcal{H}_1$ for
any $f\in\mathcal{L}_2$. Furthermore, for any $f\in\mathcal{H}_1$
%
%e11 ###
%
\begin{equation}
\int_\mathcal{T}f(t)\beta(t)\,dt=\langle Kf, \beta\rangle
_{\mathcal{H}}.\vadjust{\goodbreak}
\end{equation}
This observation allows us to prove the following result which is
important to both numerical implementation of the procedure and our
theoretical analysis.
%t1
%
\begin{theorem}
\label{th:rep}
Assume that $\ell_n$ depends on $\eta$ only through $\eta(x_1),\eta
(x_2),\ldots,\break\eta(x_n)$; then there exist $\mathbf
{d}=(d_1,\ldots,
d_{N})'\in
\mathbb{R}^N$ and $\mathbf{c}=(c_1,\ldots, c_n)'\in\mathbb{R}^n$
such that
%
%e12 ###
%
\begin{equation}
\label{eq:rep}
\hat{\beta}_{n\lambda}(t)=\sum_{k=1}^{N} d_k\xi_k(t) +\sum
_{i=1}^n c_i
(Kx_i)(t).
\end{equation}
\end{theorem}

Theorem~\ref{th:rep} is a generalization of the well-known representer
lemma for smoothing splines (Wahba, \citeyear{Wahba1990}). It
demonstrates that
although the minimization with respect to $\eta$ is taken over an
infinite-dimensional space, the solution can actually be found in a
finite-dimensional subspace, and it suffices to evaluate the
coefficients $\mathbf{c}$ and $\mathbf{d}$ in (\ref{eq:rep}). Its
proof follows a
similar argument as that of Theorem 1.3.1 in Wahba (\citeyear
{Wahba1990}) where
$\ell_n$ is assumed to be squared error, and is therefore omitted here
for brevity.

Consider, for example, the squared error loss. The regularized
estimator is given by
%
%e13 ###
%
\begin{equation}
\label{eq:pls0}
 (\hat{\alpha}_{n\lambda},\hat{\beta}_{n\lambda} )
=\mathop{\arg\min}_{\alpha\in\mathbb{R}, \beta\in\mathcal{H}}
\Biggl\{{1\over n}\sum_{i=1}^n\biggl[y_i- \biggl(\alpha+\int
_\mathcal{T}x_i(t)\beta(t)\,dt \biggr) \biggr]^2+\lambda J(\beta)
\Biggr
\}.\hspace*{-30pt}
\end{equation}
It is not hard to see that
%
%e14 ###
%
\begin{equation}
\hat{\alpha}_{n\lambda}= \bar y - \int_\mathcal{T}\bar x(t) \hat
{\beta
}_{n\lambda}(t)\,dt,
\end{equation}
where $\bar{x}(t)={1\over n} \sum_{i=1}^n x_i(t)$ and $\bar
{y}={1\over
n} \sum_{i=1}^n y_i$ are the sample average of $x$ and $y$, respectively.
Consequently, (\ref{eq:pls0}) yields
%
%e15 ###
%
\begin{equation}
\hspace*{30pt}\hat{\beta}_{n\lambda}=\mathop{\arg\min}_{\beta
\in\mathcal{H}}
\Biggl\{{1\over n}\sum_{i=1}^n
\biggl[(y_i-\bar{y})-\int_\mathcal{T}\bigl(x_i(t)-\bar{x}(t)\bigr
)\beta
(t)\,dt \biggr]^2+\lambda J(\beta) \Biggr\}.
\end{equation}
For the purpose of illustration, assume that $\mathcal{H}=\mathcal
{W}_2^2$ and
$J(\beta)=\int(\beta'')^2$. Then $\mathcal{H}_0$ is the linear
space spanned
by $\xi_1(t)=1$ and $\xi_2(t)=t$. A popular reproducing kernel
associated with $\mathcal{H}_1$ is
%
%e16 ###
%
\begin{equation}
K(s,t)={1\over(2!)^2}B_2(s)B_2(t)-{1\over4!}B_4(|s-t|),
\end{equation}
where $B_m(\cdot)$ is the $m$th Bernoulli polynomial. The readers are
referred to Wahba (\citeyear{Wahba1990}) for further details.
Following Theorem~\ref{th:rep}, it suffices to consider $\beta$ of
the following form:
%
%e17 ###
%
\begin{equation}
{\beta}(t)=d_1 + d_2 t + \sum_{i=1}^n c_i \int_{\mathcal{T}}
[x_i(s)-\bar{x}(s)]K(t,s)\,ds
\label{beta.w2}
\end{equation}
for some $\mathbf{d}\in\mathbb{R}^2$ and $\mathbf{c}\in\mathbb
{R}^n$. Correspondingly,
\begin{eqnarray*}
&&\int_\mathcal{T}[X(t)-\bar{x}(t)]\beta(t)\,dt
\\
&&\qquad=d_1 \int_\mathcal{T}[X(t)-\bar{x}(t)]\,dt + d_2\int
_\mathcal{T}
[X(t)-\bar{x}(t)]t\,dt \\
&&\qquad\quad{}+\sum_{i=1}^n c_i \int_\mathcal{T}\int_{\mathcal{T}}
[x_i(s)-\bar
{x}(s)]K(t,s)[X(t)-\bar{x}(t)]\,ds\,dt.
\end{eqnarray*}
Note also that for $\beta$ given in (\ref{beta.w2})
%
%e18 ###
%
\begin{equation}
J({\beta})=\mathbf{c}'\Sigma\mathbf{c},
\end{equation}
where $\Sigma=(\Sigma_{ij})$ is a $n\times n$ matrix with
%
%e19 ###
%
\begin{equation}
\Sigma_{ij}=\int_\mathcal{T}\int_{\mathcal{T}} [x_i(s)-\bar
{x}(s)]K(t,s)[x_j(t)-\bar
{x}(t)]\,ds\,dt.
\end{equation}
Denote by $T=(T_{ij})$ an $n\times2$ matrix whose $(i,j)$ entry is
%
%e20 ###
%
\begin{equation}
T_{ij}=\int[x_i(t)-\bar{x}(t)]t^{j-1}\,dt
\end{equation}
for $j=1,\; 2$. Set $\mathbf{y}=(y_1,\ldots,y_n)'$.
Then
%
%e21 ###
%
\begin{equation}
\label{eq:lsvec}
\ell_n(\eta)+\lambda J(\beta)={1\over n} \Vert\mathbf{y}-
(T\mathbf{d}+
\Sigma\mathbf{c}) \Vert^2_{\ell_2}+\lambda\mathbf{c}'\Sigma
\mathbf{c},
\end{equation}
which is quadratic in $\mathbf{c}$ and $\mathbf{d}$, and the explicit
form of the
solution can be easily obtained for such a problem.
This computational problem is similar to that behind the smoothing splines.
Write $W=\Sigma+n\lambda I$; then the minimizer of (\ref{eq:lsvec}) is
given by
\begin{eqnarray*}
\mathbf{d}&=& (T'W^{-1}T )^{-1}T'W^{-1}\mathbf{y},\\
\mathbf{c}&=&W^{-1} [I-T (T'W^{-1}T )^{-1}T'W^{-1} ]\mathbf{y}.
\end{eqnarray*}

%s3 ###
\section{Simultaneous diagonalization}
\label{diagonal.sec}

Before studying the asymptotic properties of the regularized
estimators $\hat{\eta}_{n\lambda}$
and $\hat\beta_{n\lambda}$, we first investigate the relationship
between the eigen structures of the covariance operator for $X(\cdot)$
and the reproducing kernel of the functional space $\mathcal{H}$. As observed
in earlier studies [e.g., Cai and Hall (\citeyear{CaiHall2006});
Hall and Horowitz (\citeyear{HallHorowitz2007})],
eigen structures play prominent roles in determining the nature of the
estimation problem in functional linear regression.

Recall that $K$ is the reproducing kernel of $\mathcal{H}_1$. Because
$K$ is
continuous and square integrable, it follows from Mercer's
theorem [Riesz and Sz-Nagy (\citeyear{RieszSznagy1955})] that $K$
admits the following
spectral decomposition:
%
%e22 ###
%
\begin{equation}
K(s,t)=\sum_{k=1}^\infty\rho_k \psi_k(s)\psi_k(t).
\end{equation}
Here $\rho_1\ge\rho_2\ge\cdots$ are the eigenvalues of $K$, and
$\{\psi_1, \psi_2, \ldots\}$ are the corresponding eigenfunctions,
that is,
%
%e23 ###
%
\begin{equation}
K\psi_k=\rho_k \psi_k, \qquad k=1,2,\ldots.
\end{equation}
Moreover,
%
%e24 ###
%
\begin{equation}
\langle\psi_i, \psi_j \rangle_{\mathcal{L}_2}=\delta_{ij}
\quad\mbox{and}\quad\langle\psi_i, \psi_j \rangle_{K}=\delta
_{ij}/\rho_j,
\end{equation}
where $\delta_{ij}$ is the Kronecker's delta.

Consider, for example, the univariate Sobolev space $\mathcal{W}_2^m([0,1])$
with norm (\ref{eq:sobnorm}) and penalty (\ref{eq:sobpen}). Observe that
%
%e25 ###
%
\begin{equation}
\mathcal{H}_1= \biggl\{f\in\mathcal{H}\dvtx\int f^{(k)}=0,
k=0,1,\ldots
,m-1\biggr\}.
\end{equation}
It is known that [see, e.g., Wahba (\citeyear{Wahba1990})]
%
%e26 ###
%
\begin{equation}
\label{eq:sobrk}
K(s,t)={1\over(m!)^2}B_m(s)B_m(t)+{(-1)^{m-1}\over(2m)!}B_{2m}(|s-t|).
\end{equation}
Recall that $B_m$ is the $m$th Bernoulli polynomial. It is known [see,
e.g., Micchelli and Wahba (\citeyear{MicchelliWahba1981})] that in
this case, $\rho_k\asymp
k^{-2m}$, where for two
positive sequences $a_k$ and $b_k$, $a_k\asymp b_k$ means that
$a_k/b_k$ is bounded away from $0$ and
$\infty$ as $k \to\infty$.

Denote by $C$ the covariance operator for $X$, that is,
%
%e27 ###
%
\begin{equation}
C(s,t)=E \{[X(s)-E(X(s))][X(t)-E(X(t))] \}.
\end{equation}
There is a duality between reproducing kernel Hilbert spaces and
covariance operators [Stein (\citeyear{Stein1999})]. Similarly to the
reproducing kernel
$K$, assuming that the covariance operator $C$ is continuous and square
integrable, we also have the following spectral decomposition
%
%e28 ###
%
\begin{equation}
C(s,t)=\sum_{k=1}^\infty\mu_k \phi_k(s)\phi_k(t),
\end{equation}
where $\mu_1\ge\mu_2\ge\cdots$ are the eigenvalues and $\{\phi
_1,\phi
_2,\ldots\}$ are the eigenfunctions such that
%
%e29 ###
%
\begin{equation}
C\phi_k:=\int_\mathcal{T}C(\cdot, t)\phi_k(t)\,dt=\mu_k\phi
_k,\qquad
k=1,2,\ldots.
\end{equation}

The decay rate of the eigenvalues $\{\mu_k\dvtx k\ge1\}$ can be
determined by the smoothness of the covariance operator $C$. More
specifically, when $C$ satisfies the so-called Sacks--Ylvisaker
conditions of order $s$ where $s$ is a nonnegative integer
[Sacks and Ylvisaker (\citeyear
{SacksYlvisaker1966,SacksYlvisaker1968,SacksYlvisaker1970})], then
$\mu_k\asymp k^{-2(s+1)}$. The
readers are referred to the original papers by Sacks and Ylvisaker or
a more recent paper by Ritter, Wasilkowski and Wo\'zniakwski (\citeyear
{RitterWasilkowskiWozniakowski1995}) for
detailed discussions of the Sacks--Ylvisaker conditions.
The conditions are also stated in the \hyperref[Appendix]{Appendix}
for \mbox{completeness}.
Roughly speaking, a
covariance operator $C$ is said to satisfy the Sacks--Ylvisaker
conditions of order $0$ if it is twice differentiable when $s\neq t$
but not differentiable when $s=t$. A~covariance operator $C$ satisfies
the Sacks--Ylvisaker conditions of order $r$ for an integer $r>0$ if
$\partial^{2r} C(s,t)/(\partial s^r\,\partial t^r)$ satisfies the
Sacks--Ylvisaker conditions of order $0$. In this paper, we say a
covariance operator $C$ satisfies the Sacks--Ylvisaker conditions if
$C$ satisfies the Sacks--Ylvisaker conditions of order $r$ for some
$r\ge0$. Various examples of covariance functions are known to satisfy
Sacks--Ylvisaker conditions. For example, the Ornstein--Uhlenbeck
covariance function $C(s,t)=\exp(-|s-t|)$ satisfies the
Sacks--Ylvisaker conditions of order $0$. Ritter, Wasilkowski and Wo\'
zniakowski (\citeyear{RitterWasilkowskiWozniakowski1995})
recently showed that covariance functions satisfying
the Sacks--Ylvisaker conditions are also intimately related to Sobolev
spaces, a fact that is useful for the purpose of simultaneously
diagonalizing $K$ and $C$ as we shall see later.

Note that the two sets of eigenfunctions $\{\psi_1, \psi_2, \ldots\}$
and $\{\phi_1,\phi_2,\ldots\}$ may differ from each other. The two
kernels $K$ and $C$ can, however, be simultaneously
diagonalized. To avoid ambiguity, we
shall assume in what follows that $Cf\neq0$ for any $f\in\mathcal
{H}_0$ and
$f\neq0$. When using the squared error loss, this is
also a necessary condition to ensure that $E\ell_n(\eta)$ is uniquely
minimized even if $\beta$ is known to come from the finite-dimensional
space $\mathcal{H}_0$. Under this assumption, we can define a norm $\|
\cdot
\|
_R$ in $\mathcal{H}$ by
%
%e30 ###
%
\begin{equation}
\label{eq:Rrk}
\|f\|_R^2=\langle Cf, f\rangle_{\mathcal{L}_2}+J(f)=\int_{\mathcal
{T}\times
\mathcal{T}
}f(s)C(s,t)f(t)\,ds\,dt+J(f).
\end{equation}
Note that $\|\cdot\|_R$ is a norm because $\|f\|_R^2$ defined above is
a quadratic form and is zero if and only if $f=0$.

The following proposition shows that when this condition holds,
$\|\cdot\|_R$ is well defined on $\mathcal{H}$ and equivalent to its
original norm, $\|\cdot\|_{\mathcal{H}}$, in that there exist constants
$0<c_1<c_2<\infty$ such that $c_1\|f\|_R\le\|f\|_\mathcal{H}\le
c_2\|f\|_R$ for all $f\in\mathcal{H}$. In particular, $\|f\|_R<\infty$
if and only if $\|f\|_\mathcal{H}<\infty$.
%p2
%
\begin{proposition}
\label{prop:Rnorm}
If $Cf\neq0$ for any $f\in\mathcal{H}_0$ and $f\neq0$, then $\|
\cdot\|_R$
and $\|\cdot\|_{\mathcal{H}}$ are equivalent.
\end{proposition}

Let $R$ be the reproducing kernel associated with $\|\cdot\|_R$. Recall
that $R$ can also be viewed as a positive operator. Denote by $\{(\rho
'_k,\psi'_k)\dvtx k\ge1\}$ the eigenvalues and eigenfunctions of
$R$. Then
$R$ is a linear map from $\mathcal{L}_2$ to $\mathcal{L}_2$ such that
%
%e31 ###
%
\begin{equation}
R\psi'_k=\int_\mathcal{T}R(\cdot, t)\psi'_k(t)\,dt=\rho'_k\psi
_k',\qquad
k=1,2,\ldots.
\end{equation}
The square root of the positive definite operator can therefore be
given as the linear map from $\mathcal{L}_2$ to $\mathcal{L}_2$ such that
%
%e32 ###
%
\begin{equation}
R^{1/2}\psi'_k=(\rho'_k)^{1/2}\psi_k',\qquad k=1,2,\ldots.\vadjust{\goodbreak}
\end{equation}
Let $\nu_1\ge
\nu_2\ge\cdots$ be the eigenvalues of the bounded linear operator
$R^{1/2}CR^{1/2}$ and $\{\zeta_k\dvtx k=1,2, \ldots\}$ be the corresponding
orthogonal eigenfunctions in $\mathcal{L}_2$. Write $\omega_k=\nu
^{-1/2}_kR^{1/2}\zeta_k$, $k=1,2,\ldots.$ Also let $\langle\cdot
,\cdot
\rangle_R$ be the inner product associated with $\|\cdot\|_R$, that is,
for any $f,g\in\mathcal{H}$,
%
%e33 ###
%
\begin{equation}
\langle f, g\rangle_R =\tfrac14 (\|f+g\|_{R}^2-\|f-g\|_{R}^2 ).
\end{equation}
It is not hard to see that
%
%e34 ###
%
\begin{equation}
\hspace*{30pt}\langle\omega_j,\omega_k\rangle_R =\nu_j^{-1/2}\nu
_k^{-1/2}\langle
R^{1/2}\zeta_j,R^{1/2}\zeta_k\rangle_R=\nu_k^{-1}\langle\zeta
_j,\zeta
_k\rangle_{\mathcal{L}_2}=\nu_k^{-1}\delta_{jk},
\end{equation}
and
\begin{eqnarray*}
\langle C^{1/2}\omega_j,C^{1/2}\omega_k\rangle_{\mathcal
{L}_2}&=&\nu
_j^{-1/2}\nu_k^{-1/2}\langle C^{1/2}R^{1/2}\zeta
_j,C^{1/2}R^{1/2}\zeta
_k\rangle_{\mathcal{L}_2}\\
&=&\nu_j^{-1/2}\nu_k^{-1/2}\langle R^{1/2}CR^{1/2}\zeta_j,\zeta
_k\rangle
_{\mathcal{L}_2}\\
&=&\delta_{jk}.
\end{eqnarray*}
The following theorem shows that quadratic forms $\|f\|_R^2=\langle
f,f\rangle_R$ and $\langle Cf,\break f\rangle_{\mathcal{L}_2}$ can be
simultaneously diagonalized on the basis of $\{\omega_k: k\ge1\}$.
%t3
%
\begin{theorem}
\label{th:simdiag}
For any $f\in\mathcal{H}$,
%
%e35 ###
%
\begin{equation}
f=\sum_{k=1}^\infty f_k\omega_k,
\end{equation}
in the absolute sense where $f_k=\nu_k\langle f, \omega_k\rangle_R$.
Furthermore, if $\gamma_k=(\nu_k^{-1}-1)^{-1}$, then
%
%e36 ###
%
\begin{equation}
\langle f,f\rangle_R = \sum_{k=1}^\infty(1+\gamma_k^{-1}
)f_k^2\quad\mbox{and}\quad\langle Cf,f\rangle_{\mathcal{L}_2} =
\sum
_{k=1}^\infty f_k^2.
\end{equation}
Consequently,
%
%e37 ###
%
\begin{equation}
J(f)=\langle f,f\rangle_R-\langle Cf,f\rangle_{\mathcal{L}_2}=\sum
_{k=1}^\infty\gamma_k^{-1}f_k^2.
\end{equation}
\end{theorem}

Note that $\{(\gamma_k,\omega_k)\dvtx k\ge1\}$ can be determined jointly
by $\{(\rho_k,\psi_k)\dvtx k\ge1\}$ and $\{(\mu_k,\phi_k)\dvtx
k\ge1\}$.
However, in general, neither $\gamma_k$ nor $\omega_k$ can be given in
explicit form of $\{(\rho_k,\psi_k)\dvtx k\ge1\}$ and
$\{(\mu_k,\phi_k)\dvtx
k\ge1\}$. One notable exception is the case when the operators $C$ and
$K$ are commutable. In particular, the setting $\psi_k=\phi_k$,
$k=1,2,\ldots,$ is commonly adopted when studying FPCA-based approaches
[see, e.g., Cai and Hall (\citeyear{CaiHall2006});
Hall and Horowitz (\citeyear{HallHorowitz2007})].
%p4
%
\begin{proposition}
\label{prop:simdiag}
Assume that $\psi_k=\phi_k$, $k=1,2,\ldots,$ then $\gamma_k=\rho
_k\mu
_k$ and $\omega_k=\mu_k^{-1/2}\psi_k$.
\end{proposition}

In general, when $\psi_k$ and $\phi_k$ differ, such a relationship no
longer holds. The following theorem reveals that similar asymptotic
behavior of $\gamma_k$ can still be expected in many practical settings.
%t5
%
\begin{theorem}
\label{th:eigen}
Consider the one-dimensional case when $\mathcal{T}=[0,1]$. If
$\mathcal{H}$ is the
Sobolev space $\mathcal{W}^m_2([0,1])$ endowed with norm (\ref{eq:sobnorm}),
and $C$ satisfies the Sacks--Ylvisaker conditions, then $\gamma
_k\asymp
\mu_k\rho_k$.
\end{theorem}

Theorem~\ref{th:eigen} shows that under fairly general conditions
$\gamma_k\asymp\mu_k\rho_k$. In this case, there is little difference
between the general situation and the special case when $K$ and $C$
share a common set of eigenfunctions when working with the system $\{
(\gamma_k,\omega_k), k=1,2,\ldots\}$. This observation is crucial for
our theoretical development in the next section.

%s4 ###
\section{Convergence rates}
\label{rate.sec}

We now turn to the asymptotic properties of the smoothness regularized
estimators. To fix ideas, in what follows, we shall focus on the
squared error loss. Recall that in this case
%
%e38 ###
%
\begin{equation}
\hspace*{25pt} (\hat{\alpha}_{n\lambda},\hat{\beta}_{n\lambda}
)=\mathop{\arg\min}
_{\alpha\in\mathbb{R}, \beta\in\mathcal{H}}
\Biggl\{{1\over n}\sum_{i=1}^n
\biggl[y_i- \biggl(\alpha+\int_\mathcal{T}x_i(t)\beta(t)\,dt
\biggr)
\biggr]^2+\lambda
J(\beta) \Biggr\}.
\end{equation}
As shown before, the slope function can be equivalently defined as
%
%e39 ###
%
\begin{equation}
\hspace*{25pt}\hat{\beta}_{n\lambda}=\mathop{\arg\min}_{\beta
\in\mathcal{H}}
\Biggl\{{1\over n}\sum
_{i=1}^n \biggl[(y_i-\bar{y})-\int_\mathcal{T}\bigl(x_i(t)-\bar
{x}(t)\bigr)\beta
(t)\,dt \biggr]^2+\lambda J(\beta) \Biggr\},
\end{equation}
and once $\hat{\beta}_{n\lambda}$ is computed, $\hat{\alpha
}_{n\lambda
}$ is given by
%
%e40 ###
%
\begin{equation}
\hat{\alpha}_{n\lambda}=\bar{y} - \int_\mathcal{T}\bar{x}(t)\hat
{\beta
}_{n\lambda}(t)\,dt.
\end{equation}
In light of this fact, we shall focus our attention on $\hat{\beta
}_{n\lambda}$ in the following discussion for brevity. We shall also
assume that the eigenvalues of the reproducing kernel $K$ satisfies
$\rho_k\asymp k^{-2r}$ for some $r>1/2$.
Let $\mathcal{F}(s,M,K)$ be the collection of the distributions $F$ of the
process $X$ that satisfy the following conditions:
\begin{enumerate}[(a)]
\item[(a)] The eigenvalues $\mu_k$ of its covariance operator
$C(\cdot,\cdot)$ satisfy $\mu_k\asymp k^{-2s}$ for some $s>1/2$.
\item[(b)] For any function $f\in\mathcal{L}_2(\mathcal{T})$,
%
%e41 ###
%
\begin{eqnarray}
&&E \biggl(\int_\mathcal{T}f(t)[X(t)-E(X)(t)]\,dt \biggr
)^4\nonumber
\\[-8pt]\\[-8pt]
&&\qquad\le M \biggl[E \biggl(\int_\mathcal{T}f(t)[X(t)-E(X)(t)]\,dt
\biggr)^2 \biggr]^2.\nonumber
\end{eqnarray}
\item[(c)] When simultaneously diagonalizing $K$ and $C$, $\gamma
_k\asymp\rho_k\mu_k$, where $\nu_k=(1+\gamma_k^{-1})^{-1}$ is the
$k$th largest eigenvalue of $R^{1/2}CR^{1/2}$ where $R$ is the
reproducing kernel associated with $\|\cdot\|_R$ defined by (\ref{eq:Rrk}).
\end{enumerate}

The first condition specifies the smoothness of the sample path of
$X(\cdot)$. The second condition concerns the fourth moment of a linear
functional of $X(\cdot)$. This condition is satisfied with $M=3$ for a
Gaussian process because $\int f(t)X(t)\,dt$ is normally distributed. In
the light of Theorem~\ref{th:eigen}, the last condition is satisfied by
any covariance function that satisfies the Sacks--Ylvisaker conditions
if $\mathcal{H}$ is taken to be $\mathcal{W}_2^m$ with norm (\ref
{eq:sobnorm}). It
is also trivially satisfied if the eigenfunctions of the covariance
operator $C$ coincide with those of $K$.

%s4.1 ###
\subsection{Optimal rates of convergence}

We are now ready to state our main results on the optimal rates of convergence,
which are given in terms of a class of intermediate norms between
$\|f\|_K$ and
%
%e42 ###
%
\begin{equation}
\biggl(\int\!\!\int f(s)C(s,t)f(t)\,ds\,dt \biggr)^{1/2},
\end{equation}
which enables a unified treatment of both
the prediction and estimation problems. For $0\le a\le1$ define the norm
$\|\cdot\|_a$ by
%
%e43 ###
%
\begin{equation}
\Vert f\Vert_a^2=\sum_{k=1}^\infty(1+\gamma_k^{-a})f_k^2,
\end{equation}
where $f_k=\nu_k \langle f, \omega_k\rangle_R$ as shown in
Theorem~\ref{th:simdiag}.
Clearly $\|f\|_0$ reduces to $\langle Cf, f\rangle_{\mathcal{L}_2}$ whereas
$\|f\|_1=\|f\|_R$. The convergence rate
results given below are valid for all $0\le a \le1$. They cover a
range of interesting cases including the prediction error and
estimation error.

The following result gives the optimal rate of convergence for the
regularized estimator $\hat{\beta}_{n\lambda}$ with an appropriately
chosen tuning parameter $\lambda$ under the loss $\|\cdot\|_a$.
%t6
%
\begin{theorem}
\label{th:main}
Assume that $E(\varepsilon_i)=0$ and $\operatorname{Var}(\varepsilon
_i)\le M_2$. Suppose
the eigenvalues $\rho_k$ of the reproducing kernel $K$ of the RKHS
$\mathcal{H}$ satisfy $\rho_k \asymp k^{-2r}$ for some $r > 1/2$.
Then the regularized estimator $\hat\beta_{n\lambda}$ with
%
%e44 ###
%
\begin{equation}
\label{eq:lamopt}
\lambda\asymp n^{-2(r+s)/(2(r+s)+1)}
\end{equation}
satisfies
%
%e45 ###
%
\begin{eqnarray}
\label{ubd.eq}
\hspace*{30pt}&&\lim_{D\to\infty} \mathop{\overline{\lim}}_{n\to\infty} \sup
_{F\in\mathcal{F}
(s,M,K),\beta_0\in
\mathcal{H}} P \bigl(\|\hat{\beta}_{n\lambda}-\beta_0\|
_a^2>Dn^{-\afrac{2(1-a)(r+s)}{2(r+s)+1}} \bigr)\nonumber
\\[-8pt]\\[-8pt]
&&\qquad=0.\nonumber
\end{eqnarray}
\end{theorem}

Note that the rate of the optimal choice of $\lambda$ does not depend
on $a$. Theorem~\ref{th:main} shows that the optimal rate of\vadjust{\goodbreak}
convergence for the
regularized estimator $\hat{\beta}_{n\lambda}$ is
$n^{-\afrac{2(1-a)(r+s)}{2(r+s)+1}}$. The following lower bound result
demonstrates that this rate of convergence is indeed optimal among all
estimators, and consequently the upper bound in equation (\ref{ubd.eq})
cannot be improved. Denote by $\mathcal{B}$ the collection of all measurable
functions of the observations $(X_1,Y_1),\ldots, (X_n,Y_n)$.
%t7
\begin{theorem}
\label{th:main1}
Under the assumptions of Theorem~\ref{th:main}, there exists a constant
$d>0$ such that
%
%e46 ###
%
\begin{equation}
\hspace*{25pt}\mathop{\underline{\lim}}_{n\to\infty}\inf_{\tilde{\beta}\in
\mathcal{B}}\sup_{F\in
\mathcal{F}
(s,M,K),\beta_0\in\mathcal{H}}
P \bigl(\|\tilde{\beta}-\beta_0\|_a^2>d n^{-\afrac{2(1-a)(r+s)}{
2(r+s)+1}} \bigr)>0.
\end{equation}
Consequently, the regularized estimator
$\hat{\beta}_{n\lambda}$ with $\lambda\asymp n^{-2(r+s)/(2(r+s)+1)}$
is rate optimal.
\end{theorem}

The results, given in terms of $\|\cdot\|_a$, provide a wide range of
measures of the quality of an estimate for $\beta_0$. Observe that
%
%e47 ###
%
\begin{equation}
\Vert\tilde{\beta}-\beta_0 \Vert_0^2=E_{X^\ast} \biggl(\int
\tilde{\beta}(t)X^\ast(t)\,dt-\int\beta_0(t)X^\ast(t)\,dt \biggr)^2,
\end{equation}
where $X^\ast$ is an independent copy of $X$, and the expectation on
the right-hand side is taken over $X^\ast$. The right-hand side is
often referred to as the prediction error in regression. It measures
the mean squared prediction error for a random future observation on
$X$. From Theorems~\ref{th:main} and~\ref{th:main1}, we have the
following corollary.
%cor8
\begin{corollary}
\label{co:me}
Under the assumptions of Theorem~\ref{th:main}, the mean squared
optimal prediction
error of a slope function estimator over $F\in\mathcal{F}(s,M,K)$ and
$\beta_0\in\mathcal{H}$ is of the order $n^{-{2(r+s)\over
2(r+s)+1}}$ and
it can be achieved by the regularized estimator $\hat{\beta
}_{n\lambda}$
with $\lambda$ satisfying (\ref{eq:lamopt}).
\end{corollary}

The result shows that the faster the eigenvalues of the covariance
operator $C$ for $X(\cdot)$ decay, the smaller the prediction error.

When $\psi_k=\phi_k$, the prediction error of a slope function
estimator $\tilde{\beta}$ can
also be understood as the squared prediction error for a fixed
predictor $x^\ast(\cdot)$ such that $|\langle x^\ast,
\phi_k\rangle_{\mathcal{L}_2}|\asymp k^{-s}$ following the discussed
from the
last section. A similar prediction
problem has also been considered by Cai and Hall (\citeyear
{CaiHall2006}) for FPCA-based
approaches. In particular, they established a similar minimax lower
bound and showed that the lower bound can be achieved by the FPCA-based
approach, but with additional assumptions that $\mu_k-\mu_{k+1}\ge
C_0^{-1}k^{-2s-1}$, and $2r>4s+3$. Our results here indicate that both
restrictions are unnecessary for establishing the minimax rate for the
prediction error. Moreover, in contrast to the FPCA-based approach, the
regularized estimator $\hat{\beta}_{n\lambda}$ can achieve the optimal
rate without the extra requirements.\vadjust{\goodbreak}

To illustrate the generality of our results, we
consider an example where $\mathcal{T}=[0,1]$, $\mathcal{H}=\mathcal
{W}_2^{m}([0,1])$
and the stochastic process $X(\cdot)$ is a Wiener process. It is not
hard to see that the covariance operator of $X$, $C(s,t)=\min\{s,t\}$,
satisfies the Sacks--Ylvisaker conditions of order $0$ and therefore
$\mu_k\asymp k^{-2}$. By Corollary~\ref{co:me}, the minimax rate of the
prediction error in estimating $\beta_0$ is $n^{-\vafrac{2m+2}{r2m+3}}$.
Note that the condition $2r>4s+3$ required by Cai and Hall (\citeyear
{CaiHall2006}) does
not hold here for $m\le7/2$.\vspace*{-2pt}

%s4.2 ###
\subsection{\texorpdfstring{The special case of $\phi_k=\psi_k$}{The special case of phi k = psi k}}

It is of interest to further look into the case when the operators $C$ and
$K$ share a common set of eigenfunctions. As discussed in the last
section, we have in this case $\phi_k=\psi_k$ and $\gamma_k\asymp
k^{-2(r+s)}$ for all $k\ge1$. In this context, Theorems~\ref{th:main}
and~\ref{th:main1} provide bounds for more general prediction problems.
Consider estimating $\int x^\ast\beta_0$ where $x^\ast$ satisfies $
|\langle x^\ast, \phi_k\rangle_{\mathcal{L}_2} |\asymp k^{-s+q}$ for
some $0<q< s-1/2$. Note that $q<s-1/2$ is needed to ensure that $x^\ast
$ is square integrable. The squared prediction error
%
%e48 ###
%
\begin{equation}
\label{eq:prederr}
\biggl(\int\tilde{\beta}(t)x^\ast(t)\,dt-\int\beta_0(t)x^\ast(t)\,dt \biggr)^2
\end{equation}
is therefore equivalent to
$\Vert\tilde{\beta}-\beta_0\Vert_{(s-q)/(r+s)}$. The
following result is a direct consequence of Theorems~\ref{th:main} and
\ref{th:main1}.\vspace*{-2pt}
%cor9
\begin{corollary}
Suppose $x^\ast$ is a function satisfying $ |\langle x^\ast, \phi
_k\rangle_{\mathcal{L}_2} |\asymp k^{-s+q}$ for some $0<q< s-1/2$. Then
under the assumptions of Theorem~\ref{th:main},
%
%e49 ###
%
\begin{eqnarray}
\hspace*{30pt}&&\mathop{\underline{\lim}}_{n\to\infty}\inf_{\tilde{\beta}\in
\mathcal{B}}
\sup_{F\in\mathcal{F}(s,M,K),\beta_0\in\mathcal{H}} P
\biggl\{ \biggl(\int\tilde{\beta}(t)x^\ast(t)\,dt-\int\beta_0(t)x^\ast(t)\,dt \biggr)^2\nonumber
\\[-9pt]\\[-9pt]
&&\hspace*{195pt}>dn^{-\afrac{2(r+q)}{2(r+s)+1}} \biggr\}>0\nonumber
\end{eqnarray}
for some constant $d>0$, and the regularized estimator $\hat{\beta
}_{n\lambda}$ with $\lambda$
satisfying (\ref{eq:lamopt}) achieves the optimal rate of convergence
under the prediction error (\ref{eq:prederr}).\vspace*{-2pt}
\end{corollary}

It is also evident that when $\psi_k=\phi_k$, $\|\cdot\|_{s/(r+s)}$ is
equivalent to $\|\cdot\|_{\mathcal{L}_2}$. Therefore, Theorems~\ref{th:main}
and~\ref{th:main1} imply the following result.\vspace*{-2pt}
%cor10
\begin{corollary}
If $\phi_k=\psi_k$ for all $k\ge1$, then under the assumptions of
Theorem~\ref{th:main}
%
%e50 ###
%
\begin{equation}
\mathop{\underline{\lim}}_{n\to\infty}\inf_{\tilde{\beta}\in
\mathcal{B}}
\sup_{F\in\mathcal{F}(s,M,K),\beta_0\in\mathcal{H}} P \bigl( \Vert
\tilde{\beta}-\beta_0 \Vert_{\mathcal{L}_2}^2>d n^{-\afrac{2r}{2(r+s)+1}} \bigr)>0
\end{equation}
for some constant $d>0$, and the regularized estimate $\hat{\beta
}_{n\lambda}$ with $\lambda$ satisfying (\ref{eq:lamopt}) achieves the
optimal rate.\vspace*{-2pt}
\end{corollary}

This result demonstrates that the faster the eigenvalues of the
covariance operator for $X(\cdot)$ decay, the larger the estimation
error. The behavior of the estimation error thus differs significantly
from that of prediction error.\vadjust{\goodbreak}

Similar results on the lower bound have recently been obtained by
Hall and Horowitz (\citeyear{HallHorowitz2007}) who considered
estimating $\beta_0$ under the
assumption that $ |\langle\beta_0,\phi_k\rangle_{\mathcal{L}_2} |$
decays in a polynomial order. Note that this slightly differs from our
setting where $\beta_0\in\mathcal{H}$ means that
%
%e51 ###
%
\begin{equation}
\label{eq:hilcond}
\sum_{k=1}^\infty\rho_k^{-1}\langle\beta_0,\psi_k\rangle
_{\mathcal{L}
_2}^2=\sum_{k=1}^\infty\rho_k^{-1}\langle\beta_0,\phi_k\rangle
_{\mathcal{L}
_2}^2<\infty.
\end{equation}
Recall that $\rho_k\asymp k^{-2r}$. Condition (\ref{eq:hilcond}) is
comparable to, and slightly stronger than,
%
%e52 ###
%
\begin{equation}
\label{eq:hhcond}
|\langle\beta_0,\phi_k\rangle_{\mathcal{L}_2} |\le M_0k^{-r-1/2}
\end{equation}
for some constant $M_0>0$. When further assuming that $2s+1<2r$, and
$\mu_k-\mu_{k+1}\ge M_0^{-1}k^{-2s-1}$ for all $k\ge1$,
Hall and Horowitz (\citeyear{HallHorowitz2007}) obtain the same lower
bound as ours. However,
we do not require that $2s+1<2r$ which in
essence states that $\beta_0$ is smoother than the sample path of
$X$. Perhaps, more\vspace*{1pt} importantly, we do not require the spacing condition
$\mu_k-\mu_{k+1}\ge M_0^{-1}k^{-2s-1}$ on the eigenvalues because we do
not need to estimate the
corresponding eigenfunctions. Such a condition is impossible to verify
even for a standard RKHS.

%s4.3 ###
\subsection{Estimating derivatives}

Theorems~\ref{th:main} and~\ref{th:main1} can also be used for
estimating the derivatives of $\beta_0$. A natural estimator of the
$q$th derivative of $\beta_0$, $\beta_0^{(q)}$, is $\hat{\beta
}_{n\lambda}^{(q)}$, the $q$th derivative of $\hat{\beta}_{n\lambda}$.
In addition to $\phi_k=\psi_k$, assume that $ \Vert\psi
_k^{(q)}/\psi
_k \Vert_\infty\asymp k^{q}$. This clearly holds when $\mathcal
{H}=\mathcal{W}
_2^m$. In this case
%
%e53 ###
%
\begin{equation}
\big\|\tilde{\beta}^{(q)}-\beta_0^{(q)} \big\|_{\mathcal{L}_2}\le C_0\|
\tilde{\beta}-\beta_0\|_{(s+q)/(r+s)}.
\end{equation}
The following is then a direct consequence of Theorems~\ref{th:main}
and~\ref{th:main1}.
%cor11
\begin{corollary}
Assume that $\phi_k=\psi_k$ and $ \Vert\psi_k^{(q)}/\psi_k
\Vert_\infty\asymp k^{q}$ for all $k\ge1$.
Then under the assumptions of Theorem~\ref{th:main}, for some constant
$d>0$,
%
%e54 ###
%
\begin{eqnarray}
\hspace*{25pt}&&\mathop{\underline{\lim}}_{n\to\infty} \inf_{\tilde{\beta
}^{(q)}\in\mathcal{B}}\sup
_{F\in\mathcal{F}
(s,M,K),\beta_0\in\mathcal{H}} P \bigl( \big\|\tilde{\beta}^{(q)}-\beta_0^{(q)}\big \|_{\mathcal{L}_2}^2
>dn^{-\afrac{2(r-q)}{2(r+s)+1}} \bigr)\nonumber
\\[-8pt]\\[-8pt]
&&\qquad>0,\nonumber
\end{eqnarray}
and the regularized estimate $\hat{\beta}_{n\lambda}$ with $\lambda$
satisfying (\ref{eq:lamopt}) achieves the optimal rate.
\end{corollary}

Finally, we note that although we have focused on the squared error loss
here, the method of regularization can be easily extended to handle
other goodness of fit measures as well as the generalized functional
linear regression [Cardot and Sarda (\citeyear{CardotSarda2005}) and
M\"uller and Stadtm\"uller (\citeyear{MullerStadtmuller2005})]. We
shall leave these extensions for future studies.

%s5 ###
\section{Numerical results}
\label{numerical.sec}

The Representer Theorem given in Section~\ref{representer.sec} makes
the regularized estimators easy to implement. Similarly to smoothness
regularized estimators in other contexts [see, e.g., Wahba (\citeyear{Wahba1990})],
$\hat{\eta}_{n\lambda}$ and $\hat{\beta}_{n\lambda}$ can be expressed
as a linear combination of a
finite number of known basis functions although the minimization in
(\ref{eq:mor}) is taken over an infinitely-dimensional space. Existing
algorithms for smoothing splines can thus be used to compute our
regularized estimators $\hat{\eta}_{n\lambda}$, $\hat{\beta
}_{n\lambda
}$ and
$\hat\alpha_{n\lambda}$.

To demonstrate the merits of the proposed estimators in finite sample
settings, we carried out a set of simulation studies. We adopt the
simulation setting of Hall and Horowitz (\citeyear{HallHorowitz2007})
where $\mathcal{T}=[0,1]$. The
true slope function $\beta_0$ is given by
%
%e55 ###
%
\begin{equation}
\beta_0=\sum_{k=1}^{50} 4(-1)^{k+1}k^{-2} \phi_k,
\end{equation}
where $\phi_1(t)=1$ and $\phi_{k+1}(t)=\sqrt{2}\cos(k\pi t)$ for
$k\ge
1$. The random function $X$ was generated as
%
%e56 ###
%
\begin{equation}
X=\sum_{k=1}^{50} \zeta_kZ_k \phi_k,
\end{equation}
where $Z_k$ are independently sampled from the uniform distribution
on\break
$[-\sqrt{3},\sqrt{3}]$ and $\zeta_k$ are deterministic. It is not hard
to see that $\zeta_k^2$ are the eigenvalues of the covariance function
of $X$. Following Hall and Horowitz (\citeyear{HallHorowitz2007}), two
sets of $\zeta_k$ were
used. In the first set, the eigenvalues are well spaced: $\zeta
_k=(-1)^{k+1}k^{-\nu/2}$ with $\nu=1.1, 1.5, 2$ or $4$. In the second set,
%
%e57 ###
%
\begin{equation}
\hspace*{20pt}\zeta_k= \cases{
1,&\quad$k=1$,\cr
0.2(-1)^{k+1}(1-0.0001k),&\quad$2\le k\le4$,\cr
0.2(-1)^{k+1} [(5\lfloor k/5\rfloor)^{-\nu/2}-0.0001(k \operatorname{mod}\ 5) ],&\quad$k\ge5$.
}
\end{equation}

As in Hall and Horowitz (\citeyear{HallHorowitz2007}), regression
models with $\varepsilon\sim
N(0,\sigma^2)$ where $\sigma=0.5$ and $1$ were considered. To comprehend
the effect of sample size, we consider $n=50, 100, 200$ and $500$. We
apply the regularization method to each simulated dataset and examine
its estimation accuracy as measured by integrated squared error $\|\hat
{\beta}_{n\lambda}-\beta_0\|_{\mathcal{L}_2}^2$ and prediction
error $\|
\hat
{\beta}_{n\lambda}-\beta_0\|_{0}^2$. For the purpose of illustration,
we take $\mathcal{H}=\mathcal{W}_2^2$ and $J(\beta)=\int(\beta'')^2$, for which
the detailed estimation procedure is given in Section~\ref{representer.sec}. For each
setting, the experiment was repeated 1000 times.

As is common in most smoothing methods, the choice of the tuning
parameter plays an important role in the performance of the
regularized estimators.
Data-driven choice of the tuning parameter is a difficult problem.
Here we apply the commonly used practical strategy of empirically
choosing the value of $\lambda$ through the generalized cross
validation. Note that the regularized estimator is a linear estimator
in that $\hat{\mathbf{y}}=H(\lambda)\mathbf{y}$ where
$\hat{\mathbf{y}}=(\hat{\eta}_{n\lambda}(x_1),\ldots, \hat{\eta
}_{n\lambda}(x_n))'$
and $H(\lambda)$ is the so-called hat matrix\vadjust{\goodbreak} depending on $\lambda$.
We then select the tuning parameter $\lambda$ that minimizes
%
%e58 ###
%
\begin{equation}
\operatorname{GCV}(\lambda)={(1/ n) \Vert\hat{\mathbf{y}}-\mathbf
{y}\Vert_{\ell
_2}^2\over(1-\operatorname{tr}(H(\lambda))/n )^2}.
\end{equation}
Denote by $\hat{\lambda}^{\mathrm{GCV}}$ the resulting choice of the tuning
parameter.

%f1 ###
%
\begin{figure}

\includegraphics{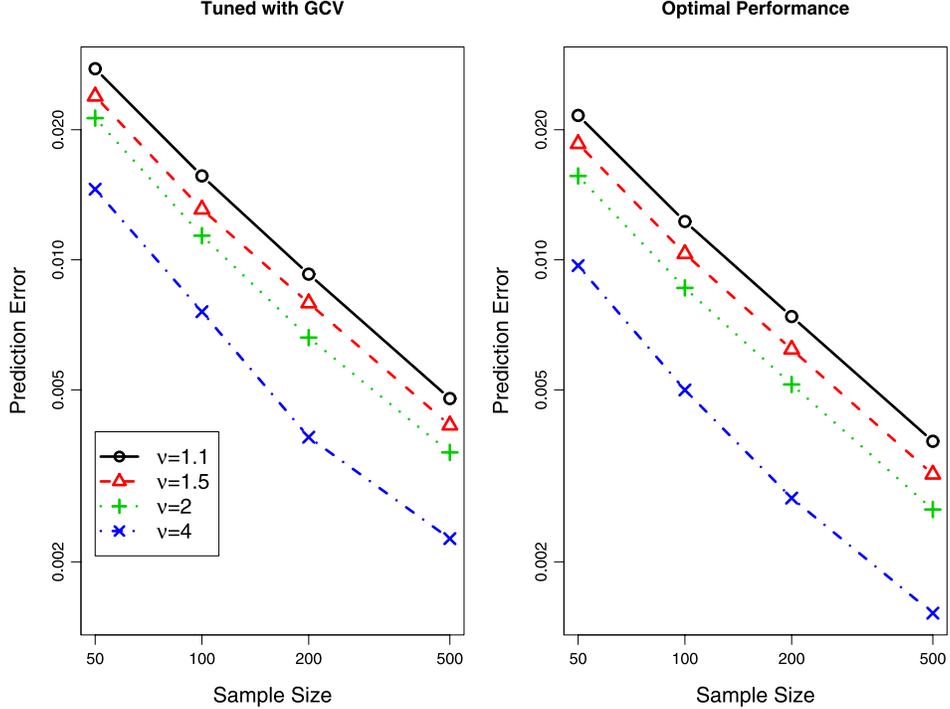}

\caption{Prediction errors of the regularized estimator ($\sigma=0.5$):
$X$ was simulated with a covariance function with well-spaced
eigenvalues. The results are averaged over 1000 runs. Black solid
lines, red dashed lines, green dotted lines and blue dash-dotted lines
correspond to $\nu= 1.1,\; 1.5,\; 2$ and $4$, respectively. Both axes
are in log scale.}
\label{fig:sim-well-pred-small}
\end{figure}

We begin with the setting of well-spaced eigenvalues. The left panel
of Figure~\ref{fig:sim-well-pred-small} shows the prediction error,
$\|\hat{\beta}_{n\lambda}-\beta_0\|_{0}^2$, for each combination
of $\nu$ value and sample size when $\sigma=0.5$. The results were
averaged over 1000 simulation runs in each setting. Both axes are
given in the log scale. The plot suggests that the estimation error
converges at a polynomial rate as sample size $n$ increases, which
agrees with our theoretical results from the previous section.
Furthermore, one can observe that with the same sample size, the
prediction error tends to be smaller for larger $\nu$. This also
confirms our theoretical development which indicates that the faster
the eigenvalues of the covariance operator for $X(\cdot)$ decay, the
smaller the prediction error.\looseness=1

%f2 ###
%
\begin{figure}

\includegraphics{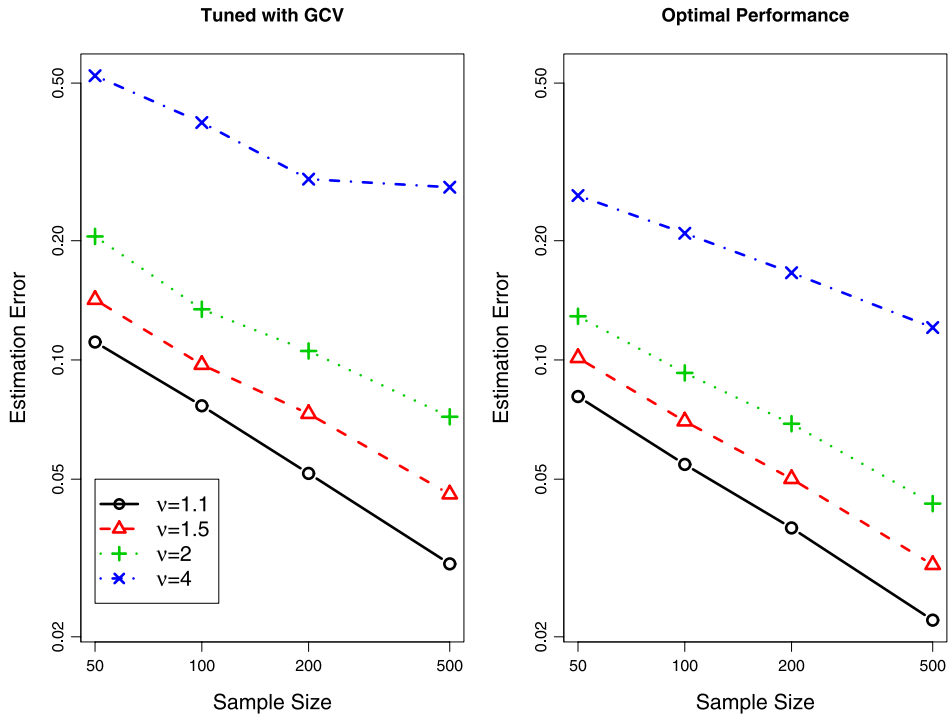}

\caption{Estimation errors of the regularized estimator ($\sigma=0.5$):
$X$ was simulated with a covariance function with well-spaced
eigenvalues. The results are averaged over 1000 runs. Black solid
lines, red dashed lines, green dotted lines and blue dash-dotted lines
correspond to $\nu= 1.1, 1.5, 2$ and $4$, respectively. Both axes are
in log scale.}
\label{fig:sim-well-est-small}\vspace*{-3pt}
\end{figure}

To better understand the performance of the smoothness regularized
estimator and the GCV choice of the tuning parameter, we also recorded
the performance of an oracle estimator whose tuning parameter is
chosen to minimize the prediction error. This choice of the tuning
parameter ensures the optimal performance of the regularized
estimator. It is, however, noteworthy that this is not a legitimate
statistical estimator since it depends on the knowledge of unknown
slope function $\beta_0$. The right panel of Figure~\ref
{fig:sim-well-pred-small} shows the prediction error associated with
this choice of tuning parameter. It behaves similarly to the estimate
with $\lambda$ chosen by GCV. Note that the comparison between the two
panels suggest that GCV generally leads to near optimal performance.

We now turn to the estimation error. Figure
\ref{fig:sim-well-est-small} shows the estimation errors, averaged
over 1000 simulation runs, with $\lambda$ chosen by GCV or minimizing
the estimation error for each combination of sample size and $\nu$
value. Similarly to the prediction error, the plots suggest a polynomial
rate of convergence of the estimation error when the sample size
increases, and GCV again leads to near-optimal choice of the tuning
parameter.\vadjust{\goodbreak}

%f3 ###
%
\begin{figure}

\includegraphics{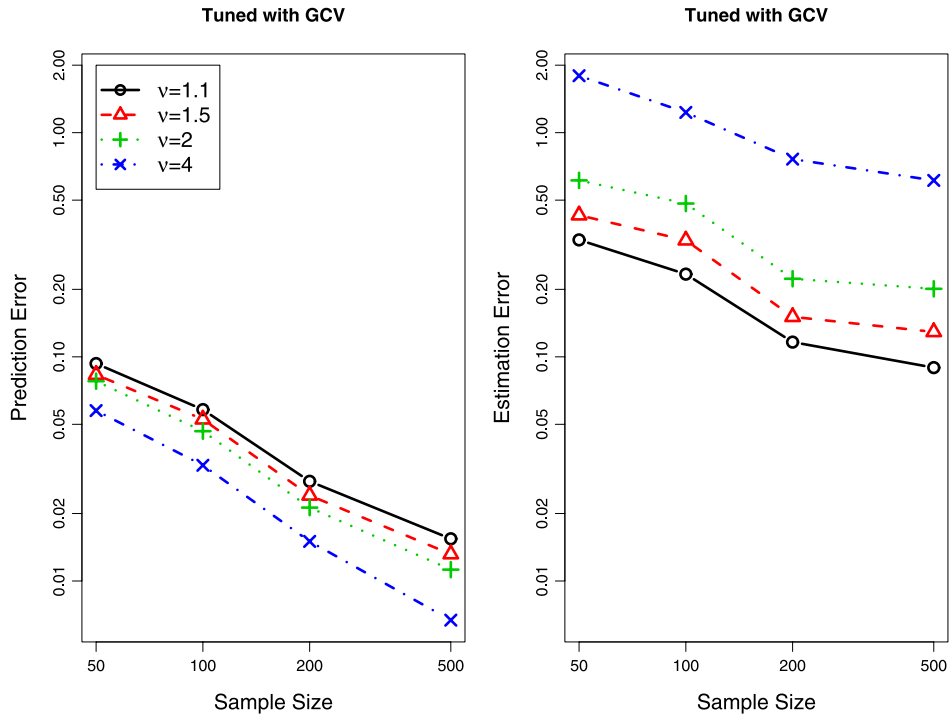}

\caption{Estimation and prediction errors of the regularized estimator
($\sigma^2=1^2$): $X$ was simulated with a covariance function with
well-spaced eigenvalues. The results are averaged over 1000 runs. Black
solid lines, red dashed lines, green dotted lines and blue dash-dotted
lines correspond to $\nu= 1.1, 1.5, 2$ and $4$, respectively. Both axes
are in log scale.}
\label{fig:sim-well-large}\vspace*{-3pt}
\end{figure}

A comparison between Figures~\ref{fig:sim-well-pred-small} and
\ref{fig:sim-well-est-small} suggests that when $X$ is smoother
(larger $\nu$), prediction (as measured by the prediction error) is
easier, but estimation (as measured by the estimation error) tends to
be harder, which highlights the difference between prediction and
estimation in functional linear regression. We also note that this
observation is in agreement with our theoretical results from the
previous section where it is shown that the estimation error decreases
at the rate of $n^{-2r/(2(r+s)+1)}$ which decelerates as $s$ increases;
whereas the prediction error decreases at the rate of
$n^{-2(r+s)/(2(r+s)+1)}$ which accelerates as $s$ increases.

Figure~\ref{fig:sim-well-large} reports the prediction and estimation
error when tuned with GCV for the large noise ($\sigma=1$) setting.
Observations similar to those for the small noise setting can also be
made. Furthermore, notice that the prediction errors are much smaller
than the estimation error, which confirms our finding from the previous
section that prediction is an easier problem in the context of
functional linear regression.

\begin{figure}

\includegraphics{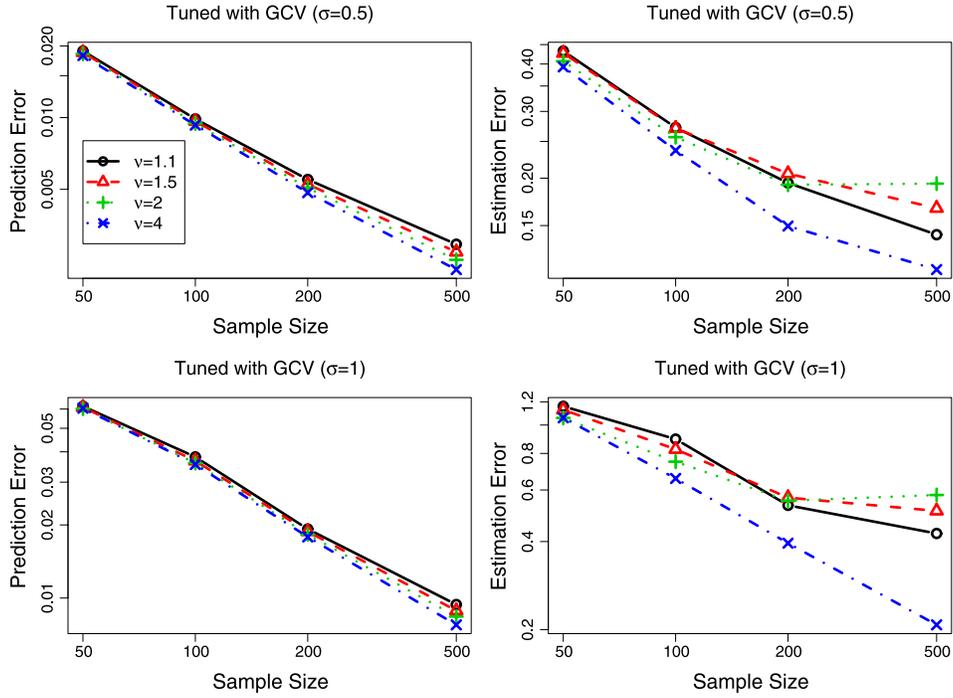}

\caption{Estimation and prediction errors of the regularized estimator:
$X$ was simulated with a covariance function with closely-spaced
eigenvalues. The results are averaged over 1000 runs. Both axes are in
log scale. Note that $y$-axes are of different scales across panels.}
\label{fig:close}
\end{figure}

The numerical results in the setting with closely spaced eigenvalues
are qualitatively similar to those\vadjust{\goodbreak} in the setting with well-spaced eigenvalues.
Figure~\ref{fig:close} summarizes the results obtained for the setting
with closely spaced eigenvalues.

%f4 ###
%

We also note that the performance of the regularization estimate with
$\lambda$ tuned with GCV compares favorably with those from Hall and
Horowitz (\citeyear{HallHorowitz2007})
using FPCA-based methods even though their results are
obtained with optimal rather than data-driven choice of the tuning parameters.

%s6 ###
\section{Proofs}\label{proof.sec}

%s6.1 ###
\subsection{\texorpdfstring{Proof of Proposition~\protect\ref{prop:Rnorm}}{Proof of Proposition 2}}

Observe that
%
%e59 ###
%
\begin{equation}
\int_{\mathcal{T}\times\mathcal{T}}f(s)C(s,t)f(t)\,ds\,dt\le\mu_1\|f\|_{\mathcal{L}_2}^2\le c_1\|f\|_\mathcal{H}^2
\end{equation}
for some constant $c_1>0$. Together with the fact that $J(f)\le\|f\|
_\mathcal{H}^2$, we conclude that
%
%e60 ###
%
\begin{equation}
\|f\|_R^2= \int_{\mathcal{T}\times\mathcal
{T}}f(s)C(s,t)f(t)\,ds\,dt+J(f)\le
(c_1+1)\|
f\|_\mathcal{H}^2.
\end{equation}

Recall that $\xi_k$, $k=1,\ldots, N,$ are the orthonormal basis of
$\mathcal{H}
_0$. Under the assumption of the proposition, the matrix $\Sigma
=(\langle C\xi_j,\xi_k\rangle_\mathcal{H})_{1\le j,k\le N}$ is a positive definite
matrix. Denote
by $\mu_1'\ge\mu_2'\ge\cdots\ge\mu_N'>0$ its eigenvalues. It is clear
that for any $f_0\in\mathcal{H}_0$
%
%e61 ###
%
\begin{equation}
\|f_0\|_R^2\ge\mu_N'\|f_0\|_\mathcal{H}^2.
\end{equation}
Note also that for any $f_1\in\mathcal{H}_1$,
%
%e62 ###
%
\begin{equation}
\|f_1\|_\mathcal{H}^2=J(f_1)\le\|f_1\|_R^2.
\end{equation}

For any $f\in\mathcal{H}$, we can write $f:=f_0+f_1$ where $f_0\in
\mathcal{H}_0$
and $f_1\in\mathcal{H}_1$. Then
%
%e63 ###
%
\begin{equation}
\|f\|_R^2=\int_{\mathcal{T}\times\mathcal{T}}f(s)C(s,t)f(t)\,ds\,dt+\|
f_1\|_{\mathcal{H}}^2.
\end{equation}
Recall that
%
%e64 ###
%
\begin{equation}
\|f_1\|_{\mathcal{H}}^2\ge\rho_1^{-1}\|f_1\|_{\mathcal{L}_2}^2\ge
\rho
_1^{-1}\mu
_1^{-1}\int_{\mathcal{T}\times\mathcal{T}}f_0(s)C(s,t)f_0(t)\,ds\,dt.
\end{equation}
For brevity, assume that $\rho_1=\mu_1=1$ without loss of generality.
By the Cauchy--Schwarz inequality,
\begin{eqnarray*}
\|f\|_R^2&\ge&{1\over2}\int_{\mathcal{T}\times\mathcal{T}
}f(s)C(s,t)f(t)\,ds\,dt+\|
f_1\|_{\mathcal{H}}^2\\
&\ge&{1\over2}\int_{\mathcal{T}\times\mathcal{T}
}f_0(s)C(s,t)f_0(t)\,ds\,dt+{3\over
2}\int_{\mathcal{T}\times\mathcal{T}}f_1(s)C(s,t)f_1(t)\,ds\,dt\\
&&{}- \biggl(\int_{\mathcal{T}\times\mathcal{T}}f_0(s)C(s,t)f_0(t)\,ds\,dt \biggr)^{1/2}
\\
&&\hspace*{12pt}{}\times\biggl(\int_{\mathcal{T}\times\mathcal{T}}f_1(s)C(s,t)f_1(t)\,ds\,dt \biggr)^{1/2}\\
&\ge&{1\over3}\int_{\mathcal{T}\times\mathcal{T}}f_0(s)C(s,t)f_0(t)\,ds\,dt,
\end{eqnarray*}
where we used the fact that $3a^2/2-ab\ge-b^2/6$ in deriving the last
inequality. Therefore,
%
%e65 ###
%
\begin{equation}
{3\over\mu_N'}\|f\|_R^2\ge\|f_0\|_\mathcal{H}^2.
\end{equation}
Together with the facts that $\|f\|_\mathcal{H}^2=\|f_0\|_\mathcal
{H}^2+\|f_1\|
_\mathcal{H}^2$ and
%
%e66 ###
%
\begin{equation}
\|f\|_R^2\ge J(f_1)\ge\|f_1\|_\mathcal{H}^2,
\end{equation}
we conclude that
%
%e67 ###
%
\begin{equation}
\|f\|_R^2\ge(1+3/\mu_N')^{-1}\|f\|_\mathcal{H}^2.
\end{equation}
The proof is now complete.

%s6.2 ###
\subsection{\texorpdfstring{Proof of Theorem~\protect\ref{th:simdiag}}{Proof of Theorem 3}}
First note that
\begin{eqnarray*}
R^{-1/2}f&=&\sum_{k=1}^\infty\langle R^{-1/2}f, \zeta_k\rangle
_{\mathcal{L}
_2} \zeta_k
=\sum_{k=1}^\infty\langle R^{-1/2}f, \nu_k^{1/2}R^{-1/2}\omega
_k\rangle
_{\mathcal{L}_2} \nu_k^{1/2}R^{-1/2}\omega_k\\
&=&R^{-1/2} \Biggl(\sum_{k=1}^\infty\nu_k\langle R^{-1/2}f,
R^{-1/2}\omega_k\rangle_{\mathcal{L}_2} \omega_k \Biggr)
=R^{-1/2} \Biggl(\sum_{k=1}^\infty\nu_k\langle f, \omega_k\rangle_R
\omega_k \Biggr).
\end{eqnarray*}
Applying bounded positive definite operator $R^{1/2}$ to both sides
leads to
%
%e68 ###
%
\begin{equation}
f=\sum_{k=1}^\infty\nu_k\langle f, \omega_k\rangle_R \omega_k.
\end{equation}

Recall that $\langle\omega_k, \omega_j\rangle_R=\nu_k^{-1}\delta
_{kj}$. Therefore,
\begin{eqnarray*}
\|f\|_R^2&=& \Bigg\langle\sum_{k=1}^\infty\nu_k\langle f, \omega
_k\rangle_R \omega_k, \sum_{j=1}^\infty\nu_j\langle f, \omega
_j\rangle
_R \omega_j \Bigg\rangle_R
\\
&=&\sum_{k,j=1} \nu_k\nu_j\langle f, \omega_k\rangle_R\langle f,
\omega
_j\rangle_R\langle\omega_k,\omega_j\rangle_R\\
&=&\sum_{k=1} \nu_k \langle f, \omega_k\rangle_R^2.
\end{eqnarray*}

Similarly, because $\langle C\omega_k, \omega_j\rangle_{\mathcal{L}
_2}=\delta_{kj}$,
\begin{eqnarray*}
\langle Cf, f\rangle_{\mathcal{L}_2}
&=& \Bigg\langle C \Biggl(\sum_{k=1}^\infty
\nu_k\langle f, \omega_k\rangle_R \omega_k \Biggr), \sum_{j=1}^\infty
\nu
_j\langle f, \omega_j\rangle_R \omega_j \Bigg\rangle_{\mathcal{L}_2}\\
&=& \Bigg\langle\sum_{k=1}^\infty\nu_k\langle f, \omega_k\rangle_R
C\omega_k, \sum_{j=1}^\infty\nu_j\langle f, \omega_j\rangle_R
\omega
_j \Bigg\rangle_{\mathcal{L}_2}\\
&=&\sum_{k,j=1} \nu_k\nu_j\langle f, \omega_k\rangle_R\langle f,
\omega
_j\rangle_R\langle C\omega_k,\omega_j\rangle_{\mathcal{L}_2}\\
&=&\sum_{k=1} \nu_k^2 \langle f, \omega_k\rangle_R^2.
\end{eqnarray*}

%s6.3 ###
\subsection{\texorpdfstring{Proof of Proposition~\protect\ref{prop:simdiag}}{Proof of Proposition 4}}

Recall that for any $f\in\mathcal{H}_0$, $Cf\neq0$ if and only if $f=0$,
which implies that $\mathcal{H}_0\cap{\rm l.s.}\{\phi_k\dvtx  k\ge1\}
^{\perp
}=\{
0\}$. Together with the fact that $\mathcal{H}_0\cap\mathcal{H}_1=\{
0\}$, we
conclude that $\mathcal{H}=\mathcal{H}_1={\rm l.s.}\{\phi_k\dvtx  k\ge1\}
$. It is not
hard to see that for any $f,g\in\mathcal{H}$,
%
%e69 ###
%
\begin{equation}
\langle f,g\rangle_R=\int_{\mathcal{T}\times\mathcal{T}
}f(s)C(s,t)g(t)\,ds\,dt+\langle
f, g\rangle_K.
\end{equation}
In particular,
%
%e70 ###
%
\begin{equation}
\langle\psi_j,\psi_k\rangle_R= (\mu_k+\rho_k^{-1} )\delta_{jk},
\end{equation}
which implies that $\{((\mu_k+\rho_k^{-1})^{-1},\psi_k)\dvtx  k\ge1\}$ is
also the eigen system of $R$, that is,
%
%e71 ###
%
\begin{equation}
R(s,t)=\sum_{k=1}^\infty(\mu_k+\rho_k^{-1})^{-1}\psi_k(s)\psi_k(t).
\end{equation}
Then
%
%e72 ###
%
\begin{equation}
R\psi_k:=\int_\mathcal{T}R(\cdot, t)\psi_k(t)\,dt=(\mu_k+\rho
_k^{-1})^{-1}\psi
_k,\qquad k=1,2,\ldots.
\end{equation}
Therefore,
\begin{eqnarray*}
R^{1/2}CR^{1/2}\psi_k&=&R^{1/2}C \bigl((\mu_k+\rho_k^{-1})^{-1/2}\psi
_k \bigr)\\
&=&R^{1/2} \bigl(\mu_k(\mu_k+\rho_k^{-1})^{-1/2}\psi_k \bigr)
= (1+\rho_k^{-1}\mu_k^{-1} )^{-1} \psi_k,
\end{eqnarray*}
which implies that $\zeta_k=\psi_k=\phi_k$, $\nu_k= (1+\rho
_k^{-1}\mu_k^{-1} )^{-1}$ and $\gamma_k=\rho_k\mu_k$. Consequently,
%
%e73 ###
%
\begin{equation}
\omega_k=\nu_k^{-1/2}R^{1/2}\psi_k=\nu_k^{-1/2}(\mu_k+\rho
_k^{-1})^{-1/2}\psi_k=\mu_k^{-1/2}\psi_k.
\end{equation}

%s6.4 ###
\subsection{\texorpdfstring{Proof of Theorem~\protect\ref{th:eigen}}{Proof of Theorem 5}}
Recall that $\mathcal{H}=\mathcal{W}_2^m$, which implies that $\rho
_k\asymp
k^{-2m}$. By Corollary 2 of Ritter, Wasilkowski and Wo\'zniakowski
(\citeyear{RitterWasilkowskiWozniakowski1995}), $\mu_k\asymp
k^{-2(s+1)}$. It therefore suffices to show $\gamma_k\asymp
k^{-2(s+1+m)}$. The key idea of the proof is a result from Ritter,
Wasilkowski and Wo\'zniakowski (\citeyear{RitterWasilkowskiWozniakowski1995})
indicating that the reproducing kernel Hilbert space
associated with $C$ differs from $\mathcal{W}^{s+1}_2([0,1])$ only by a
finite-dimensional linear space of polynomials.

Denote by $Q_r$ the reproducing kernel for $\mathcal{W}^{r}_2([0,1])$.
Observe that $Q_r^{1/2}(\mathcal{L}_2)=\mathcal{W}^r_2$ [e.g., Cucker
and Smale (\citeyear{CuckerSmale2001})]. We begin by quantifying the decay rate of
$\lambda_{k}
(Q_{m}^{1/2}Q_{s+1}Q_m^{1/2} )$. By Sobolev's embedding theorem,
$(Q_{s+1}^{1/2}Q_{m}^{1/2})(\mathcal{L}_2)=Q_{s+1}^{1/2}(\mathcal
{W}_2^{m})=\mathcal{W}
_2^{m+s+1}$. Therefore, $Q_{m}^{1/2}Q_{s+1}Q_m^{1/2}$ is equivalent to
$Q_{m+s+1}$. Denote by $\lambda_k(Q)$ be the $k$th largest eigenvalue
of a positive definite operator $Q$. Let $\{h_k\dvtx  k\ge1\}$ be the
eigenfunctions of $Q_{m+s+1}$, that is, $Q_{m+s+1}h_k=\lambda
_k(Q_{m+s+1})h_k$, $k=1,2,\ldots.$ Denote by $\mathcal{F}_k$ and
$\mathcal{F}
_k^{\perp}$ the linear space spanned by $\{h_j\dvtx  1\le j\le k\}$ and $\{
h_j\dvtx  j\ge k+1\}$, respectively. By the Courant--Fischer--Weyl min--max principle,
\begin{eqnarray*}
\lambda_k (Q_{m}^{1/2}Q_{s+1}Q_m^{1/2} )&\ge& \min_{f\in\mathcal{F}
_k} \Vert Q_{s+1}^{1/2}Q_m^{1/2}f \Vert_{\mathcal{L}_2}^2/\|f\|
_{\mathcal{L}_2}^2\\
&\ge& C_1 \min_{f\in\mathcal{F}_k} \Vert Q_{m+s+1}^{1/2}f \Vert
_{\mathcal{L}_2}^2/\|f\|_{\mathcal{L}_2}^2\\
&\ge& C_1\lambda_k(Q_{m+s+1})
\end{eqnarray*}
for some constant $C_1>0$. On the other hand,
\begin{eqnarray*}
\lambda_{k} (Q_{m}^{1/2}Q_{s+1}Q_m^{1/2} )&\le& \max_{f\in\mathcal{F}
_{k-1}^{\perp}} \Vert Q_{s+1}^{1/2}Q_m^{1/2}f \Vert_{\mathcal{L}
_2}^2/\|f\|_{\mathcal{L}_2}^2\\
&\le& C_2 \min_{f\in\mathcal{F}_{k-1}^{\perp}} \Vert
Q_{m+s+1}^{1/2}f \Vert_{\mathcal{L}_2}^2/\|f\|_{\mathcal{L}_2}^2\\
&\le& C_2\lambda_k(Q_{m+s+1})
\end{eqnarray*}
for some constant $C_2>0$. In summary, we have $\lambda_{k}
(Q_{m}^{1/2}Q_{s+1}Q_m^{1/2} )\asymp\break\times k^{-2(m+s+1)}$.

As shown by Ritter, Wasilkowski and Wo\'zniakowski [(\citeyear{RitterWasilkowskiWozniakowski1995}), Theorem 1, page~525], there exist $D$
and $U$ such that $Q_{s+1}=D+U$, $D$ has at most $2(s+1)$ nonzero
eigenvalues and $\|U^{1/2}f\|_{\mathcal{L}_2}$ is equivalent to $\|
C^{1/2}f\|
_{\mathcal{L}_2}$. Moreover, the eigenfunctions of $D$, denoted by $\{
g_1,\ldots, g_d\}$ ($d\le2(s+1)$) are polynomials of order no greater
than $2s+1$. Denote $\mathcal{G}$ the space spanned by $\{g_1,\ldots
,g_d\}$.
Clearly $\mathcal{G}\subset\mathcal{W}_2^{m}=Q_m^{1/2}(\mathcal
{L}_2)$. Denote $\{
\tilde
{h}_j \dvtx  j\ge1\}$ the eigenfunctions of $Q_{m}^{1/2}Q_{s+1}Q_m^{1/2}$.
Let $\tilde{\mathcal{F}}_k$ and $\tilde{\mathcal{F}}_k^\perp$ be
defined similarly
as $\mathcal{F}_k$ and $\mathcal{F}_k^{\perp}$. Then by the
Courant--Fischer--Weyl
min--max principle,
\begin{eqnarray*}
\lambda_{k-d} (Q_{m}^{1/2}UQ_m^{1/2} )&\ge& \min_{f\in\tilde
{\mathcal{F}
}_k\cap Q_m^{-1/2}(\mathcal{G})^{\perp}} \Vert U^{1/2}Q_m^{1/2}f
\Vert
_{\mathcal{L}
_2}^2/\|f\|_{\mathcal{L}_2}^2\\
&=& \min_{f\in\tilde{\mathcal{F}}_k\cap Q_m^{-1/2}(\mathcal
{G})^{\perp}}
\Vert
Q_{s+1}^{1/2}Q_m^{1/2}f \Vert_{\mathcal{L}_2}^2/\|f\|_{\mathcal
{L}_2}^2\\
&=& \min_{f\in\tilde{\mathcal{F}}_k} \Vert Q_{s+1}^{1/2}Q_m^{1/2}f
\Vert
_{\mathcal{L}_2}^2/\|f\|_{\mathcal{L}_2}^2\\
&\ge& C_1\lambda_k(Q_{m+s+1})
\end{eqnarray*}
for some constant $C_1>0$. On the other hand,
\begin{eqnarray*}
\lambda_{k+d} (Q_{m}^{1/2}Q_{s+1}Q_m^{1/2} )&\le& \max_{f\in\tilde
{\mathcal{F}
}_{k-1}^{\perp}\cap Q_m^{-1/2}(\mathcal{G})^{\perp}} \Vert U^{1/2}Q_m^{1/2}f
\Vert_{\mathcal{L}_2}^2/\|f\|_{\mathcal{L}_2}^2\\
&=& \max_{f\in\tilde{\mathcal{F}}_{k-1}^{\perp}\cap
Q_m^{-1/2}(\mathcal{G}
)^{\perp
}} \Vert Q_{s+1}^{1/2}Q_m^{1/2}f \Vert_{\mathcal{L}_2}^2/\|f\|
_{\mathcal{L}
_2}^2\\
&=& \min_{f\in\tilde{\mathcal{F}}_{k-1}^{\perp}} \Vert
Q_{s+1}^{1/2}Q_m^{1/2}f \Vert_{\mathcal{L}_2}^2/\|f\|_{\mathcal
{L}_2}^2\\
&\le& C_2\lambda_k(Q_{m+s+1})
\end{eqnarray*}
for some constant $C_2>0$. Hence $\lambda_{k} (Q_{m}^{1/2}UQ_m^{1/2}
)\asymp k^{-2(m+s+1)}$.

Because $Q_{m}^{1/2}UQ_m^{1/2}$ is equivalent to $R^{1/2}CR^{1/2}$,
following a similar argument as before, by the Courant--Fischer--Weyl
min--max principle, we complete the the proof.

%s6.5 ###
\subsection{\texorpdfstring{Proof of Theorem~\protect\ref{th:main}}{Proof of Theorem 6}}
We now proceed to prove Theorem~\ref{th:main}. The analysis follows a
similar spirit as the technique commonly used in the study of the rate
of convergence of smoothing splines [see, e.g., Silverman (\citeyear{Silverman1982});
Cox and O'Sullivan (\citeyear{CoxOSullivan1990})]. For brevity, we
shall assume that $EX(\cdot)=0$ in
the rest of the proof. In this case, $\alpha_0$ can be estimated by
$\bar{y}$ and $\beta_0$ by
%
%e74 ###
%
\begin{equation}
\hat{\beta}_{n\lambda}=\mathop{\arg\min}_{\beta\in\mathcal{H}}
\Biggl[{1\over n}\sum_{i=1}^n
\biggl(y_i-\int_\mathcal{T}x_i(t)\beta(t)\,dt \biggr)^2+\lambda J(\beta) \Biggr].
\end{equation}
The proof below also applies to the more general setting when $EX(\cdot
)\neq0$ but with considerable technical obscurity.

Recall that
%
%e75 ###
%
\begin{equation}
\ell_n(\beta)={1\over n}\sum_{i=1}^n \biggl(y_i-\int_\mathcal
{T}x_i(t)\beta
(t)\,dt \biggr)^2.
\end{equation}
Observe that
\begin{eqnarray*}
\ell_\infty(\beta)&:=& E\ell_n(\beta)
=E \biggl[Y-\int_\mathcal{T}X(t)\beta(t)\,dt\biggr]^2\\
&=&\sigma^2+\int_\mathcal{T}\int_\mathcal{T}[\beta(s)-\beta_0(s) ]C(s,t)[\beta(t)-\beta_0(t) ]\,ds\,dt\\
&=&\sigma^2+ \Vert\beta-\beta_0 \Vert^2_0.
\end{eqnarray*}
Write
%
%e76 ###
%
\begin{equation}
\bar{\beta}_{\infty\lambda}=\mathop{\arg\min}_{\beta\in
\mathcal{H}} \{\ell
_\infty(\beta
)+\lambda J(\beta) \}.
\end{equation}
Clearly
%
%e77 ###
%
\begin{equation}
\hat{\beta}_{n\lambda}-\beta_0= (\hat{\beta}_{n\lambda}-\bar
{\beta
}_{\infty\lambda} )+ (\bar{\beta}_{\infty\lambda}-\beta_0 ).
\end{equation}
We refer to the two terms on the right-hand side stochastic error and
deterministic error, respectively.

%s6.5.1 ###
\subsubsection{Deterministic error}
Write $\beta_0(\cdot)=\sum_{k=1}^\infty a_k\omega_k(\cdot)$ and
${\beta}(\cdot)=\break\sum_{k=1}^\infty b_k\omega_k(\cdot)$.
Then Theorem~\ref{th:simdiag} implies that
\begin{eqnarray*}
\ell_\infty({\beta})=\sigma^2+\sum_{k=1}^\infty(b_{k}-{a}_{k})^2,\qquad
J({\beta})=\sum_{k=1}^\infty\gamma_k^{-1} b_k^2.
\end{eqnarray*}
Therefore,
%
%e78 ###
%
\begin{equation}
\bar{\beta}_{\infty\lambda}(\cdot)=\sum_{k=1}^\infty{ a_k\over
1+\lambda\gamma_k^{-1}} \omega_k(\cdot)=: \sum_{k=1}^\infty\bar
{b}_k\omega_k(\cdot).
\end{equation}

It can then be computed that for any $a<1$,
\begin{eqnarray*}
\Vert\bar{\beta}_{\infty\lambda}-\beta_0 \Vert_a^2&=&\sum
_{k=1}^\infty
(1+\gamma_k^{-a})(\bar{b}_k-a_k)^2\\[-2pt]
&=&\sum_{k=1}^\infty(1+\gamma_k^{-a}) \biggl({ \lambda\gamma_k^{-1}\over
1+\lambda\gamma_k^{-1}} \biggr)^2a_k^2\\[-2pt]
&\le&\lambda^2 \sup_k {(1+\gamma^{-a})\gamma_k^{-1}\over
(1+\lambda
\gamma_k^{-1} )^2}\sum_{k=1}^\infty\gamma_k^{-1}a_k^2\\[-2pt]
&=&\lambda^2J(\beta_0) \sup_k {(1+\gamma^{-a})\gamma_k^{-1}\over
(1+\lambda\gamma_k^{-1} )^2}.
\end{eqnarray*}
Now note that
\begin{eqnarray*}
\sup_k {(1+\gamma^{-a})\gamma_k^{-1}\over(1+\lambda\gamma_k^{-1}
)^2}&\le& \sup_{x> 0} {(1+x^{-a})x^{-1}\over(1+\lambda x^{-1} )^2}\\[-2pt]
&\le&\sup_{x> 0} {x^{-1}\over(1+\lambda x^{-1} )^2}+\sup_{x> 0}
{x^{-(a+1)}\over(1+\lambda x^{-1} )^2}\\[-2pt]
&=& {1\over\inf_{x> 0} (x^{1/2}+\lambda x^{-1/2} )^2}+ {1\over\inf
_{x> 0} (x^{(a+1)/2}+\lambda x^{-(1-a)/2} )^2}\\[-2pt]
&=&{1\over4\lambda}+C_0\lambda^{-(a+1)}.
\end{eqnarray*}
Hereafter, we use $C_0$ to denote a generic positive constant. In
summary, we have
%l12
\begin{lemma}
\label{le:detererr}
If $\lambda$ is bounded from above, then
\[
\Vert\bar{\beta}_{\infty\lambda}-\beta_0 \Vert^2_{a}=O ({\lambda^{1-a}}
J(\beta_0) ).
\]
\end{lemma}

%s6.5.2 ###
\subsubsection{Stochastic error}
Next, we consider the stochastic error $\hat{\beta}_{n\lambda}-\bar
{\beta}_{\infty\lambda}$. Denote
\begin{eqnarray*}
D\ell_{n} (\beta)f
&=&-{2\over n}\sum_{i=1}^n \biggl[ \biggl(y_i-\int_\mathcal{T}x_i(t)\beta(t)\,dt \biggr)\int_\mathcal{T}x_i(t)f(t)\,dt \biggr],\\[-2pt]
D\ell_{\infty} (\beta)f&=&-2E_X \biggl(\int_\mathcal{T}X(t) [\beta_0(t)-\beta(t)]\,dt\int_\mathcal{T}X(t)f(t)\,dt \biggr)\\[-2pt]
&=&-2\int_\mathcal{T}\int_\mathcal{T}[\beta_0(s)-\beta(s)]C(s,t)f(t)\,ds\,dt,\\[-2pt]
D^2\ell_{n} (\beta)fg&=&{2\over n}\sum_{i=1}^n \biggl[\int_\mathcal{T}x_i(t)f(t)\,dt\int_\mathcal{T}x_i(t)g(t)\,dt \biggr],\\[-2pt]
D^2\ell_{\infty} (\beta)fg&=&2\int_\mathcal{T}\int_\mathcal{T}f(s)C(s,t)g(t)\,ds\,dt.
\end{eqnarray*}
Also write $\ell_{n\lambda}(\beta)=\ell_n(\beta)+\lambda J(\beta
)$ and
$\ell_{\infty\lambda}=\ell_\infty(\beta)+\lambda J(\beta)$. Denote
$G_\lambda=(1/2)D^2\ell_{\infty\lambda}(\bar{\beta}_{\infty
\lambda})$ and
%
%e79 ###
%
\begin{equation}
\tilde{\beta}=\bar{\beta}_{\infty\lambda}-\tfrac12 G_\lambda
^{-1}D\ell
_{n\lambda}(\bar{\beta}_{\infty\lambda}).
\end{equation}
It is clear that
%
%e80 ###
%
\begin{equation}
\hat{\beta}_{n\lambda}-\bar{\beta}_{\infty\lambda}= (\hat{\beta
}_{n\lambda}-\tilde{\beta} )+ (\tilde{\beta}-\bar{\beta}_{\infty
\lambda
} ).
\end{equation}
We now study the two terms on the right-hand side separately. For
brevity, we shall abbreviate the subscripts of $\hat{\beta}$ and
$\bar
{\beta}$ in what follows. We begin with $\tilde{\beta}-\bar{\beta}$.
Hereafter we shall omit the subscript for brevity if no confusion occurs.
%l13
\begin{lemma}
\label{le:stocherr1}
For any $0\le a\le1$,
%
%e81 ###
%
\begin{equation}
E \Vert\tilde{\beta}-\bar{\beta} \Vert_a^2\asymp n^{-1}\lambda^{-
(a+\afrac{1}{2(r+s)} )}.
\end{equation}
\end{lemma}

\begin{pf}
Notice that $D\ell_{n\lambda}(\bar{\beta
})=D\ell_{n\lambda}(\bar{\beta})-D\ell_{\infty\lambda}(\bar
{\beta
})=D\ell_{n}(\bar{\beta})-D\ell_{\infty}(\bar{\beta})$. Therefore
\begin{eqnarray*}
E [D\ell_{n\lambda}(\bar{\beta})f ]^2
&=&E [D\ell_{n}(\bar{\beta})f-D\ell_{\infty}(\bar{\beta})f ]^2\\[-2pt]
&=&{4\over n} \operatorname{Var}\biggl[ \biggl(Y-\int_\mathcal{T}X(t)\bar{\beta}(t)\,dt \biggr)\int_\mathcal{T}X(t)f(t)\,dt \biggr]\\[-2pt]
&\le&{4\over n} E \biggl[ \biggl(Y-\int_\mathcal{T}X(t)\bar{\beta}(t)\,dt \biggr)\int_\mathcal{T}X(t)f(t)\,dt \biggr]^2\\[-2pt]
&=&{4\over n} E \biggl(\int_\mathcal{T}X(t) [\beta_0(t)-\bar{\beta}(t)]\,dt\int_\mathcal{T}X(t)f(t)\,dt \biggr)^2\\[-2pt]
&&{}+{4\sigma^2\over n}E \biggl(\int_\mathcal{T}X(t)f(t)\,dt \biggr)^2,
\end{eqnarray*}
where we used the fact that $\varepsilon=Y-\int X\beta_0$ is
uncorrelated with $X$. To bound the first term, an application of the
Cauchy--Schwarz inequality yields
\begin{eqnarray*}
&&E \biggl(\int_\mathcal{T}X(t) [\beta_0(t)-\bar{\beta}(t) ]\,dt\int_\mathcal{T}X(t)f(t)\,dt \biggr)^2\\[-2pt]
&&\qquad\le \biggl\{E \biggl(\int_\mathcal{T}X(t) [\beta_0(t)-\bar{\beta}(t) ]\,dt \biggr)^4E\biggl(\int_\mathcal{T}X(t)f(t)\,dt \biggr)^4 \biggr\}^{1/2}\\[-2pt]
&&\qquad\le M \|\beta_0-\bar{\beta}\|^2_0\|f\|_0^2,
\end{eqnarray*}
where the second inequality holds by the second condition of $\mathcal{F}
(s,M,K)$. Therefore,
%
%e82 ###
%
\begin{equation}
E [D\ell_{n\lambda}(\bar{\beta})f ]^2\le{4M\over n} \|\beta
_0-\bar{\beta
}\|^2_0\|f\|_0^2+{4\sigma^2\over n}\|f\|_0^2,\vadjust{\goodbreak}
\end{equation}
which by Lemma~\ref{le:detererr} is further bounded by $(C_0\sigma^2/
n) \|f\|_0^2$ for some positive constant $C_0$. Recall that $\|\omega
_k\|_0=1$. We have
%
%e83 ###
%
\begin{equation}
E [D\ell_{n\lambda}(\bar{\beta})\omega_k ]^2\le C_0\sigma^2/ n.
\end{equation}
Thus, by the definition of $\tilde{\beta}$,
\begin{eqnarray*}
E \Vert\tilde{\beta}-\bar{\beta} \Vert_a^2
&=&E \bigg\Vert{1\over 2}G_\lambda^{-1}D\ell_{n\lambda}(\bar{\beta}) \bigg\Vert_a^2\\
&=& {1\over4}E \Biggl[\sum_{k=1}^\infty(1+\gamma_k^{-a})(1+\lambda\gamma_k^{-1})^{-2} (D\ell_{n\lambda}(\bar{\beta})\omega_k )^2 \Biggr]\\
&\le&{C_0\sigma^2\over4n}\sum_{k=1}^\infty(1+\gamma_k^{-a})(1+\lambda\gamma_k^{-1})^{-2}\\
&\le&{C_0\sigma^2\over4n} \sum_{k=1}^\infty\bigl(1+k^{2a(r+s)}\bigr)\bigl(1+\lambda k^{2(r+s)}\bigr)^{-2}\\
&\asymp&{C_0\sigma^2\over4n}\int_1^\infty x^{2a(r+s)}\bigl(1+\lambda x^{2(r+s)}\bigr)^{-2}\,dx\\
&\asymp&{C_0\sigma^2\over4n}\int_1^\infty\bigl(1+\lambda x^{2(r+s)/(2a(r+s)+1)} \bigr)^{-2}\,dx\\
&=&{C_0\sigma^2\over4n}\lambda^{- (a+\afrac{1}{2(r+s)} )} \int_{\lambda^{a+\afrac{1}{2(r+s)}}}^\infty\bigl(1+x^{2(r+s)/(2a(r+s)+1)} \bigr)^{-2}\,dx\\
&\asymp&n^{-1}\lambda^{- (a+\afrac{1}{2(r+s)} )}.
\end{eqnarray*}
The proof is now complete.
\end{pf}

Now we are in position to bound $E\|\hat{\beta}-\tilde{\beta}\|_a^2$.
By definition,
%
%e84 ###
%
\begin{equation}
G_\lambda(\hat{\beta}-\tilde{\beta})=\tfrac12D^2\ell_{\infty
\lambda
}(\bar{\beta})(\hat{\beta}-\tilde{\beta}).
\end{equation}
First-order condition implies that
%
%e85 ###
%
\begin{equation}
D\ell_{n\lambda}(\hat{\beta})=D\ell_{n\lambda}(\bar{\beta
})+D^2\ell
_{n\lambda}(\bar{\beta})(\hat{\beta}-\bar{\beta})=0,
\end{equation}
where we used the fact that $\ell_{n,\lambda}$ is quadratic. Together
with the fact that
%
%e86 ###
%
\begin{equation}
D\ell_{n\lambda}(\bar{\beta})+D^2\ell_{\infty\lambda}(\bar{\beta})(\tilde{\beta}-\bar{\beta})=0,
\end{equation}
we have
\begin{eqnarray*}
D^2\ell_{\infty\lambda}(\bar{\beta})(\hat{\beta}-\tilde{\beta})
&=&D^2\ell_{\infty\lambda}(\bar{\beta})(\hat{\beta}-\bar{\beta})+D^2\ell_{\infty\lambda}(\bar{\beta})(\bar{\beta}-\tilde{\beta})\\
&=&D^2\ell_{\infty\lambda}(\bar{\beta})(\hat{\beta}-\bar{\beta})-D^2\ell_{n\lambda}(\bar{\beta})(\hat{\beta}-\bar{\beta})\\
&=&D^2\ell_{\infty}(\bar{\beta})(\hat{\beta}-\bar{\beta})-D^2\ell_{n}(\bar{\beta})(\hat{\beta}-\bar{\beta}).
\end{eqnarray*}
Therefore,
%
%e87 ###
%
\begin{equation}
(\hat{\beta}-\tilde{\beta})=\tfrac12G_\lambda^{-1} [D^2\ell_{\infty}(\bar{\beta})(\hat{\beta}-\bar{\beta})-D^2\ell_{n}(\bar{\beta})(\hat{\beta}-\bar{\beta}) ].
\end{equation}
Write
%
%e88 ###
%
\begin{equation}
\hat{\beta}=\sum_{k=1}^\infty\hat{b}_k\omega_k\quad \mbox{and}\quad
\bar{\beta}=\sum_{k=0}^\infty\bar{b}_k\omega_k.
\end{equation}
Then
\begin{eqnarray*}
&&\Vert\hat{\beta}-\tilde{\beta} \Vert^2_a
\\
&&\qquad={1\over4}\sum_{k=1}^\infty(1+\lambda\gamma_k^{-1})^{-2}(1+\gamma_k^{-a})
\\
&&\qquad\quad\hspace*{23pt}{} \times\Biggl[\sum_{j=1}^\infty(\hat{b}_j-\bar{b}_j)\int_\mathcal{T}\int_\mathcal{T}
\omega_j(s) \Biggl({1\over n}\sum_{i=1}^nx_i(t)x_i(s)-C(s,t) \Biggr)
\\
&&\qquad\quad\hspace*{217pt}{}\times\omega_k(t)\,ds\,dt\Biggr]^2\\
&&\qquad\le{1\over4}\sum_{k=1}^\infty(1+\lambda\gamma_k^{-1})^{-2}(1+\gamma_k^{-a}) \Biggl[ \sum_{j=1}^\infty(\hat{b}_j-\bar{b}_j)^2(1+\gamma_j^{-c})\Biggr]\\
&&\qquad\quad\hspace*{23pt}{} \times\Biggl(\sum_{j=1}^\infty(1+\gamma_j^{-c})^{-1}
\Biggl[\int_\mathcal{T}\int_\mathcal{T}\omega_j(s) \Biggl({1\over n}\sum_{i=1}^nx_i(t)x_i(s)-C(s,t) \Biggr)
\\
&&\qquad\quad\hspace*{236pt}{}\times\omega_k(t)\,ds\,dt\Biggr]^2 \Biggr),
\end{eqnarray*}
where the inequality is due to the Cauchy--Schwarz inequality.

Note that
\begin{eqnarray*}
\hspace*{-5pt}&&E \Biggl(\sum_{j=1}^\infty(1+\gamma_j^{-c})^{-1} \Biggl[\int_\mathcal{T}\omega_j(s)\Biggl({1\over n}\sum_{i=1}^nx_i(t)x_i(s)-C(s,t) \Biggr)\omega_k(t)\,ds\,dt \Biggr]^2 \Biggr)\\
\hspace*{-5pt}&&\qquad={1\over n}\sum_{j=1}^\infty(1+\gamma_j^{-c})^{-1}\operatorname{Var}\biggl(\int_\mathcal{T}\omega_j(t)X(t)\,dt\,\int_\mathcal{T}\omega_k(t)X(t)\,dt \biggr)\\
\hspace*{-5pt}&&\qquad\le{1\over n}\sum_{j=1}^\infty(1+\gamma_j^{-c})^{-1} E \biggl[ \biggl(\int_\mathcal{T}\omega_j(t)X(t)\,dt \biggr)^2 \biggl(\int_\mathcal{T}\omega_k(t)X(t)\,dt \biggr)^2 \biggr]\\
\hspace*{-5pt}&&\qquad\le{1\over n}\sum_{j=1}^\infty(1+\gamma_j^{-c})^{-1} E \biggl[ \biggl(\int_\mathcal{T}\omega_j(t)X(t)\,dt \biggr)^4\biggr]^{1/2}
%&&\qquad\quad\hspace*{20pt}{}\times
E \biggl[ \biggl(\int_\mathcal{T}\omega_k(t)X(t)\,dt \biggr)^4\biggr]^{1/2}\\
\hspace*{-5pt}&&\qquad\le{M\over n}\sum_{j=1}^\infty(1+\gamma_j^{-c})^{-1} E \biggl[ \biggl(\int_\mathcal{T}\omega_j(t)X(t)\,dt \biggr)^2 \biggr]E \biggl[ \biggl(\int_\mathcal{T}\omega_k(t)X(t)\,dt \biggr)^2\biggr]\\
\hspace*{-5pt}&&\qquad={M\over n}\sum_{j=1}^\infty(1+\gamma_j^{-c})^{-1} \asymp n^{-1},
\end{eqnarray*}
provided that $c>1/2(r+s)$. On the other hand,
\begin{eqnarray*}
&&\sum_{k=1}^\infty(1+\lambda\gamma_k^{-1})^{-2}(1+\gamma_k^{-a})
\\
&&\qquad\le C_0\sum_{k=1}^\infty\bigl(1+\lambda k^{2(r+s)}\bigr)^{-2}\bigl(1+k^{2a(r+s)}\bigr)\\
&&\qquad\asymp\int_1^\infty\bigl(1+\lambda x^{2(r+s)} \bigr)^{-2}x^{2a(r+s)}\,dx\\
&&\qquad\asymp\int_1^\infty\bigl(1+\lambda x^{2(r+s)/(2a(r+s)+1)} \bigr)^{-2}\,dx\\
&&\qquad=\lambda^{- (a+1/(2(r+s)) )}\int_{\lambda^{a+1/(2(r+s))}}^\infty\bigl(1+x^{2(r+s)/(2a(r+s)+1)} \bigr)^{-2}\,dx\\
&&\qquad\asymp\lambda^{- (a+1/(2(r+s)) )}.
\end{eqnarray*}

To sum up,
%
%e89 ###
%
\begin{equation}
\label{eq:hat-til}
\Vert\hat{\beta}-\tilde{\beta} \Vert^2_a=O_p \bigl(n^{-1}\lambda^{-
(a+1/(2(r+s)) )} \Vert\hat{\beta}-\bar{\beta} \Vert^2_c \bigr).
\end{equation}
In particular, taking $a=c$ yields
%
%e90 ###
%
\begin{equation}
\Vert\hat{\beta}-\tilde{\beta} \Vert^2_c=O_p \bigl(n^{-1}\lambda^{-
(c+1/(2(r+s)) )} \Vert\hat{\beta}-\bar{\beta} \Vert^2_c \bigr).
\end{equation}
If
%
%e91 ###
%
\begin{equation}
n^{-1}\lambda^{- (c+1/(2(r+s)) )} \to0,
\end{equation}
then
%
%e92 ###
%
\begin{equation}
\Vert\hat{\beta}-\tilde{\beta} \Vert_c=o_p ( \Vert\hat{\beta
}-\bar
{\beta} \Vert_c ).
\end{equation}
Together with the triangular inequality
%
%e93 ###
%
\begin{equation}
\Vert\tilde{\beta}-\bar{\beta} \Vert_c \ge\Vert\hat{\beta
}-\bar
{\beta} \Vert_c- \Vert\hat{\beta}-\tilde{\beta} \Vert_c=\bigl(1-o_p(1)\bigr)
\Vert\hat{\beta}-\bar{\beta} \Vert_c.
\end{equation}
Therefore,
%
%e94 ###
%
\begin{equation}
\Vert\hat{\beta}-\bar{\beta} \Vert_c=O_p ( \Vert\tilde{\beta
}-\bar
{\beta} \Vert_c )
\end{equation}
Together with Lemma~\ref{le:stocherr1}, we have
%
%e95 ###
%
\begin{equation}
\Vert\hat{\beta}-\bar{\beta} \Vert^2_c=O_p \bigl(n^{-1}\lambda^{-
(c+\afrac{1}{2(r+s)} )} \bigr)=o_p(1).
\end{equation}

Putting it back to (\ref{eq:hat-til}), we now have:
%l14
\begin{lemma}
\label{le:stocherr2}
If there also exists some ${1/2(r+s)}<c\le1$ such that $n^{-1}\times\lambda
^{- (c+{1/2(r+s)} )}\to0$, then
%
%e96 ###
%
\begin{equation}
\Vert\hat{\beta}-\tilde{\beta} \Vert^2_a=o_p \bigl(n^{-1}\lambda^{-
(a+{1/2(r+s)} )} \bigr).
\end{equation}
\end{lemma}

Combining Lemmas~\ref{le:detererr}--\ref
{le:stocherr2}, we have
%
%e97 ###
%
\begin{eqnarray}
\hspace*{30pt}&&\lim_{D\to\infty} \mathop{\overline{\lim}}_{n\to\infty} \sup
_{F\in\mathcal{F}
(s,M,K),\beta_0\in
\mathcal{H}} P \bigl(\|\hat{\beta}_{n\lambda}-\beta_0\|
_a^2>Dn^{-\afrac{2(1-a)(r+s)}{
2(r+s)+1}} \bigr)\nonumber
\\[-8pt]\\[-8pt]
&&\qquad=0\nonumber
\end{eqnarray}
by taking $\lambda\asymp n^{-2(r+s)/(2(r+s)+1)}$.

%s6.6 ###
\subsection{\texorpdfstring{Proof of Theorem~\protect\ref{th:main1}}{Proof of Theorem 7}}

We now set out to show that $n^{-2(1-a)(r+s)/(2(r+s)+1)}$ is the
optimal rate. It follows from a similar argument as that of Hall and
Horowitz (\citeyear{HallHorowitz2007}). Consider a setting where $\psi
_k=\phi_k$, $k=1,2,\ldots
.$ Clearly in this case we also have $\omega_k=\mu_k^{-1/2}\phi_k$. It
suffices to show that the rate is optimal in this special case. Recall
that $\beta_0=\sum a_k\phi_k$. Set
%
%e98 ###
%
\begin{equation}
a_k= \cases{L_n^{-1/2} k^{-r}\theta_k, &\quad $L_{n}+1\le k\le2L_n$,\cr
0,&\quad \mbox{otherwise},
}
\end{equation}
where $L_n$ is the integer part of $n^{1/(2(r+s)+1)}$, and $\theta_k$
is either $0$ or $1$. It is clear that
%
%e99 ###
%
\begin{equation}
\|\beta_0\|_K^2\le\sum_{k=L_n+1}^{2L_n} L_n^{-1} =1.
\end{equation}
Therefore $\beta_0\in\mathcal{H}$. Now let $X$ admit the following
expansion: $X=\sum_k \xi_k k^{-s}\phi_k$ where $\xi_k$s are independent
random variables drawn from a uniform distribution on $[-\sqrt
{3},\sqrt
{3}]$. Simple algebraic manipulation shows that the distribution of $X$
belongs to $\mathcal{F}(s, 3)$. The observed data are
%
%e100 ###
%
\begin{equation}
y_i=\sum_{k=L_n+1}^{2L_n} L_n^{-1/2}k^{-(r+s)}\xi_{ik}\theta
_k+\varepsilon_i,\qquad i=1,\ldots, n,
\end{equation}
where the noise $\varepsilon_i$ is assumed to be independently sampled
from $N(0,M_2)$. As shown in Hall and Horowitz
(\citeyear{HallHorowitz2007}),
%
%e101 ###
%
\begin{equation}
\lim_{n\to\infty} \inf_{L_n<j\le2L_n}\inf_{\tilde{\theta}_j}\sup^\ast E(\tilde{\theta}_j-\theta_j)^2>0,
\end{equation}
where $\mathop{\sup^{}}\limits^{*}$
denotes the supremum over all $2^{L_n}$
choices of $(\theta_{L_n+1},\ldots, \theta_{2L_n})$, and $\inf
_{\tilde
{\theta}}$ is taken over all measurable functions $\tilde{\theta}_j$ of
the data. Therefore, for any estimate $\tilde{\beta}$,
%
%e102 ###
%
\begin{eqnarray}
\label{eq:mseldb}
\sup^\ast\|\tilde{\beta}-\beta_0\|_a^2
&=&\sup^\ast\sum_{k=L_n+1}^{2L_n}
L_n^{-1}k^{-2(1-a)(r+s)}E(\tilde{\theta}_j-\theta_j)^2\nonumber
\\[-8pt]\\[-8pt]
&\ge& Mn^{-\afrac{2(1-a)(r+s)}{2(r+s)+1}}\nonumber
\end{eqnarray}
for some constant $M>0$.

Denote
%
%e103 ###
%
\begin{equation}
\tilde{\hspace*{-2pt}\tilde{\theta}}_k=
\cases{1, &\quad $\tilde{\theta}_k>1$,\cr
\tilde{\theta}_k,&\quad $0\le\tilde{\theta}_k\le1$,\cr
0, &\quad $\tilde{\theta}_k<0$.
}
\end{equation}
It is easy to see that
%
%e104 ###
%
\begin{eqnarray}
&&\sum_{k=L_n+1}^{2L_n} L_n^{-1}k^{-2(1-a)(r+s)}(\tilde{\theta}_j-\theta_j)^2\nonumber
\\[-8pt]\\[-8pt]
&&\qquad\ge\sum_{k=L_n+1}^{2L_n} L_n^{-1}k^{-2(1-a)(r+s)}(\hspace*{2pt}\tilde{\hspace*{-2pt}\tilde{\theta}}_j-\theta_j)^2.\nonumber
\end{eqnarray}
Hence, we can assume that $0\le\tilde{\theta}_j\le1$ without loss of
generality in establishing the lower bound. Subsequently,
%
%e105 ###
%
\begin{eqnarray}
\sum_{k=L_n+1}^{2L_n} L_n^{-1}k^{-2(1-a)(r+s)}(\tilde{\theta}_j-\theta_j)^2
&\le&\sum_{k=L_n+1}^{2L_n} L_n^{-1}k^{-2(1-a)(r+s)}\nonumber
\\[-8pt]\\[-8pt]
&\le& L_n^{-2(1-a)(r+s)}.\nonumber
\end{eqnarray}
Together with (\ref{eq:mseldb}), this implies that
%
%e106 ###
%
\begin{equation}
\lim_{n\to\infty}\inf_{\tilde{\beta}}\sup^\ast
P \bigl(\Vert\tilde{\beta}-\beta\Vert_a^2>dn^{-\afrac{2(1-a)(r+s)}{2(r+s)+1}} \bigr)>0
\end{equation}
for some constant $d>0$.

\begin{appendix}
\section*{Appendix: Sacks--Ylvisaker conditions}\label{Appendix}

In Section~\ref{diagonal.sec}, we discussed the relationship between
the smoothness of $C$ and the decay of its eigenvalues. More precisely,
the smoothness can be quantified by the so-called Sacks--Ylvisaker
conditions. Following Ritter, Wasilkowski and Wo\'zniakowski (\citeyear
{RitterWasilkowskiWozniakowski1995}), denote
%
%e107 ###
%
\begin{eqnarray}
\Omega_+&=&\{(s,t)\in(0,1)^2\dvtx  s>t\}\quad \mbox{and}\nonumber
\\[-8pt]\\[-8pt]
\Omega_-&=&\{(s,t)\in(0,1)^2\dvtx  s<t\}.\nonumber
\end{eqnarray}
Let $\operatorname{cl}(A)$ be the closure of a set $A$. Suppose that $L$ is a
continuous function on $\Omega_+\cup\Omega_-$ such that $L|_{\Omega_j}$
is continuously extendable to $\operatorname{cl}(\Omega_j)$ for $j\in\{+,-\}$.
By $L_j$ we denote the extension of $L$ to $[0,1]^2$, which is
continuous on $\operatorname{cl}(\Omega_j)$, and on $[0,1]^2 \setminus\operatorname{cl}(\Omega_j)$. Furthermore write $M^{(k,l)}(s,t)=(\partial
^{k+l}/(\partial s^k \,\partial t^l))M(s,t)$. We say that a covariance
function $M$ on $[0,1]^2$ satisfies the Sacks--Ylvisaker conditions of
order $r$ if the following three conditions hold:
\begin{enumerate}[(A)]
\item[(A)] $L=M^{(r,r)}$ is continuous on $[0,1]^2$, and its partial
derivatives up to order 2 are continuous on $\Omega_+\cup\Omega_-$, and
they are continuously extendable to $\operatorname{cl}(\Omega_+)$ and $\operatorname{cl}(\Omega_-)$.
\item[(B)]
%
%e108 ###
%
\begin{equation}
\min_{0\le s\le1} \bigl(L_-^{(1,0)}(s,s)-L_+^{(1,0)}(s,s) \bigr)>0.
\end{equation}
\item[(C)] $L_+^{(2,0)}(s,\cdot)$ belongs to the reproducing kernel
Hilbert space spanned by $L$ and furthermore
%
%e109 ###
%
\begin{equation}
\sup_{0\le s\le1} \big\Vert L_+^{(2,0)}(s,\cdot) \big\Vert_L<\infty.
\end{equation}
\end{enumerate}
\end{appendix}

%suskaldyti doi

\printaddresses


\begin{thebibliography}{99}

%b1 ###
\bibitem[\protect\citeauthoryear{}{1975}]{Adams1975}
\textsc{Adams, R. A.} (1975).
\textit{Sobolev Spaces}.
Academic Press, New York.
\MR{0450957}

%b2 ###
\bibitem[\protect\citeauthoryear{}{2006}]{CaiHall2006}
\textsc{Cai, T.} and \textsc{Hall, P.} (2006).
Prediction in functional linear regression.
\textit{Ann. Statist.} \textbf{34} 2159--2179.
\MR{2291496}

%b3 ###
\bibitem[\protect\citeauthoryear{}{2003}]{CardotFerratySarda2003}
\textsc{Cardot, H., Ferraty, F.} and \textsc{Sarda, P.} (2003).
Spline estimators for the functional linear model.
\textit{Statist. Sinica} \textbf{13} 571--591.
\MR{1997162}

%b4 ###
\bibitem[\protect\citeauthoryear{}{2005}]{CardotSarda2005}
\textsc{Cardot, H.} and \textsc{Sarda, P.} (2005).
Estimation in generalized linear models for functional data via
penalized likelihood.
\textit{J.~Multivariate Anal.} \textbf{92} 24--41.
\MR{2102242}

%b5 ###
\bibitem[\protect\citeauthoryear{}{1990}]{CoxOSullivan1990}
\textsc{Cox, D. D.} and \textsc{O'Sullivan, F.} (1990).
Asymptotic analysis of penalized likelihood and related estimators.
\textit{Ann. Statist.} \textbf{18} 1676--1695.
\MR{1074429}

%b6 ###
\bibitem[\protect\citeauthoryear{}{2009}]{CrambesKneipSarda2009}
\textsc{Crambes, C., Kneip, A.} and \textsc{Sarda, P.} (2009).
Smoothing splines estimators for functional linear regression.
\textit{Ann. Statist.} \textbf{37} 35--72.
\MR{2488344}

%b7 ###
\bibitem[\protect\citeauthoryear{}{2001}]{CuckerSmale2001}
\textsc{Cucker, F.} and \textsc{Smale, S.} (2001).
On the mathematical foundations of learning.
\textit{Bull. Amer. Math. Soc.} \textbf{39} 1--49.
\MR{1864085}

%b8 ###
\bibitem[\protect\citeauthoryear{}{2006}]{FerratyVieu2006}
\textsc{Ferraty, F.} and \textsc{Vieu, P.} (2006).
\textit{Nonparametric Functional Data Analysis: Methods, Theory,
Applications and Implementations}.
Springer, New York.
\MR{2229687}

%b9 ###
\bibitem[\protect\citeauthoryear{}{2007}]{HallHorowitz2007}
\textsc{Hall, P.} and \textsc{Horowitz, J. L.} (2007).
Methodology and convergence rates for functional linear regression.
\textit{Ann. Statist.} \textbf{35} 70--91.
\MR{2332269}

%b10 ###
\bibitem[\protect\citeauthoryear{}{2002}]{James2002}
\textsc{James, G.} (2002).
Generalized linear models with functional predictors.
\textit{J.~Roy. Statist. Soc. Ser.~B} \textbf{64} 411--432.
\MR{1924298}

%b11 ###
\bibitem[\protect\citeauthoryear{}{2009}]{Johanness2009}
\textsc{Johannes, J.} (2009).
Nonparametric estimation in functional linear models with second order
stationary regressors.
Unpublished manuscript.
Available at
\href{http://arxiv.org/abs/0901.4266v1}{http://arxiv.org/}
\href{http://arxiv.org/abs/0901.4266v1}{abs/0901.4266v1}.

%b12 ###
\bibitem[\protect\citeauthoryear{}{2007}]{LiHsing2007}
\textsc{Li, Y.} and \textsc{Hsing, T.} (2007).
On the rates of convergence in functional linear regression.
\textit{J.~Multivariate Anal.} \textbf{98} 1782--1804.
\MR{2392433}

%b13 ###
\bibitem[\protect\citeauthoryear{}{1981}]{MicchelliWahba1981}
\textsc{Micchelli, C.} and \textsc{Wahba, G.} (1981).
Design problems for optimal surface interpolation. In
\textit{Approximation Theory and Applications} (Z. Ziegler, ed.) 329--347.
Academic Press, New York.
\MR{0615422}

%b14 ###
\bibitem[\protect\citeauthoryear{}{2005}]{MullerStadtmuller2005}
\textsc{M\"{u}ller, H. G.} and \textsc{Stadtm\"{u}ller, U.} (2005).
Generalized functional linear models.
\textit{Ann. Statist.} \textbf{33} 774--805.
\MR{2163159}

%b15 ###
\bibitem[\protect\citeauthoryear{}{2005}]{RamsaySilverman2005}
\textsc{Ramsay, J. O.} and \textsc{Silverman, B. W.} (2005).
\textit{Functional Data Analysis}, 2nd ed.
Springer, New York.
\MR{2168993}

%b16 ###
\bibitem[\protect\citeauthoryear{}{1955}]{RieszSznagy1955}
\textsc{Riesz, F.} and \textsc{Sz-Nagy, B.} (1955).
\textit{Functional Analysis}.
Ungar, New York.
\MR{0071727}

%b17 ###
\bibitem[\protect\citeauthoryear{}{1995}]{RitterWasilkowskiWozniakowski1995}
\textsc{Ritter, K., Wasilkowski, G.} and \textsc{Wo\'zniakowski, H.} (1995).
Multivariate integeration and approximation for random fields
satisfying Sacks--Ylvisaker conditions.
\textit{Ann. Appl. Probab.} \textbf{5} 518--540.
\MR{1336881}

%b18 ###
\bibitem[\protect\citeauthoryear{}{1966}]{SacksYlvisaker1966}
\textsc{Sacks, J.} and \textsc{Ylvisaker, D.} (1966).
Designs for regression problems with correlated errors.
\textit{Ann. Math. Statist.} \textbf{37} 66--89.
\MR{0192601}

%b19 ###
\bibitem[\protect\citeauthoryear{}{1968}]{SacksYlvisaker1968}
\textsc{Sacks, J.} and \textsc{Ylvisaker, D.} (1968).
Designs for regression problems with correlated errors; many parameters.
\textit{Ann. Math. Statist.} \textbf{39} 49--69.
\MR{0220424}

%b20 ###
\bibitem[\protect\citeauthoryear{}{1970}]{SacksYlvisaker1970}
\textsc{Sacks, J.} and \textsc{Ylvisaker, D.} (1970).
Designs for regression problems with correlated errors III.
\textit{Ann. Math. Statist.} \textbf{41} 2057--2074.
\MR{0270530}

%b21 ###
\bibitem[\protect\citeauthoryear{}{1982}]{Silverman1982}
\textsc{Silverman, B. W.} (1982).
On the estimation of a probability density function by the maximum
penalized likelihood method.
\textit{Ann. Statist.} \textbf{10} 795--810.
\MR{0663433}

%b22 ###
\bibitem[\protect\citeauthoryear{}{1999}]{Stein1999}
\textsc{Stein, M.} (1999).
\textit{Statistical Interpolation of Spatial Data: Some Theory for Kriging}.
Springer, New York.
\MR{1697409}

%b23 ###
\bibitem[\protect\citeauthoryear{}{1990}]{Wahba1990}
\textsc{Wahba, G.} (1990).
\textit{Spline Models for Observational Data}.
SIAM, Philadelphia.
\MR{1045442}

%b24 ###
\bibitem[\protect\citeauthoryear{}{2005}]{YaoMullerWang2005}
\textsc{Yao, F., M\"{u}ller, H. G.} and \textsc{Wang, J. L.} (2005).
Functional linear regression analysis for longitudinal data.
\textit{Ann. Statist.} \textbf{33} 2873--2903.
\MR{2253106}

%regularization estimators, \textit{Annals of Statistics} \textbf{16}
%694--712.
%Springer, New York.
%covariance kernels and application to second order processes,
%convergence for nonparametric regression, \textit{Annals of
%Statistics} \textbf{10} 1040--1053.
%Spaces}, Springer, New York.
%for Eigenvalue Approximation}, SIAM, Philadelphia.

\end{thebibliography}
\end{document}